# Phase-coupled Oscillators with Plastic Coupling: Synchronization and Stability

Andrey Gushchin, Enrique Mallada, and Ao Tang

**Abstract**—In this article we study synchronization of systems of homogeneous phase-coupled oscillators with plastic coupling strengths and arbitrary underlying topology. The dynamics of the coupling strength between two oscillators is governed by the phase difference between these oscillators. We show that, under mild assumptions, such systems are gradient systems, and always achieve frequency synchronization. Furthermore, we provide sufficient stability and instability conditions that are based on results from algebraic graph theory. For a special case when underlying topology is a tree, we formulate a criterion (necessary and sufficient condition) of stability of equilibria. For both, tree and arbitrary topologies, we provide sufficient conditions for phase-locking, i.e. convergence to a stable equilibrium almost surely. We additionally find conditions when the system possesses a unique stable equilibrium, and thus, almost global stability follows. Several examples are used to demonstrate variety of equilibria the system has, their dependence on system's parameters, and to illustrate differences in behavior of systems with constant and plastic coupling strengths.

## 1 INTRODUCTION

SYNCHRONIZATION of phase-coupled oscillators is an extensive topic of research that finds applications in a variety of disciplines including neuroscience [25], [34], [43], [44], physics [34], [5], mathematics [16] and engineering [10], [32]. The dynamic behavior of these systems can be quite rich. For example, the intrinsic symmetry of the network can produce multiple limit cycles or equilibria with relatively fixed phases (phase-locked) [3], and the heterogeneity in the natural oscillation frequency can lead to incoherence [11] or even chaos [38].

One of the most important properties of a system of phase-coupled oscillators is how coupling (or interaction) between oscillators is defined. The Kuramoto model [26] – a canonical model for studying synchronization phenomena – uses a trigonometric $\sin()$ *coupling function* that depends on the phase difference of the two interacting oscillators. However, broader classes of coupling functions have been considered – specially in applications to biological systems – and have proven to lead to richer varieties of dynamic behaviors [30], [20], [6], [31], [13].

Besides the coupling function, there are two additional elements that also affect the systems behavior: the coupling strength (the gain that multiplies the coupling function) and the interconnection topology (that describes who affects whom). There is a vast body of literature devoted to understanding the effect of these elements, including studies of networks with complete graph [36], graph of diameter two [15] or arbitrary topology [16], [14], [9], [11]. In general, larger (positive) coupling strength and more connected topologies tend to promote synchronization and lead to tightly grouped phase-locked solutions [14], [11]. However, when negative coupling strengths are allowed among oscillators, new stable phase-locked solutions and even non-synchronizing traveling waves can appear [19].

Interestingly, a common feature of all these studies is that the coupling strength is assumed to be fixed. Having constant coupling strengths generally simplifies the analysis and allows the theory to provide profound insights on the behavior of these systems. However, considering varying or plastic coupling strengths is more suitable for studying oscillations in neuroscience, since synaptic neural connections undergo modifications due to learning or forgetting processes [20], [41], [37]. This have motivated some recent empirical studies [8], [35], [39] that seek to understand the effect of dynamic coupling strength. However, with the exception of a few studies that consider plastic coupling strength for complete graph topologies and sinusoidal coupling, there has not been a systematic study of the dynamical properties of *plastic phase-coupled oscillators*.

The goal of this work is to develop a general analytical framework for studying systems of phase-coupled homogeneous oscillators with non-constant coupling and arbitrary underlying topology. We show by providing a Lyapunov function, that under mild conditions these systems always achieve frequency synchronization, and derive two sufficient conditions: one for showing stability, and another one for showing instability of an equilibrium. Moreover, these conditions characterize all equilibria when underlying topology is a tree graph. We further characterize the relationship between the system parameters and its behavior, as well as the range of admissible asymptotic coupling strengths. In particular, we show that for almost all choices of these parameters, the system converges to a stable equilibrium almost surely.

The structure of the article is the following. Section 2

- A. Gushchin is with the Center for Applied Mathematics, Cornell University, Ithaca, NY, 14850.
  E-mail: avg36@cornell.edu
- E. Mallada is with the Electrical and Computer Engineering Department, Johns Hopkins University, Baltimore, MD, 21218
- A. Tang is with the School of Electrical and Computer Engineering, Cornell University, Ithaca, NY, 14850.

formally describes the model, introduces necessary notation (Subsection 2.1), discusses related work and summarizes our results (Subsection 2.2). Section 3 provides several examples that motivate our study. Section 4 contains general theoretical results: in Subsection 4.1 a Lyapunov function is introduced and frequency synchronization of oscillators is shown. Then, in Subsection 4.2 we formulate stability and instability sufficient conditions. We apply these conditions in Subsection 4.3 to analyze stability of the in-phase and anti-phase equilibria. The results presented in Section 4 can be strengthen when the underlying network topology is a tree graph as shown in Section 5. In particular, we prove convergence to a stable equilibrium almost surely in the case of a general coupling (Subsection 5.1), and almost global stability in the case of strictly attractive or repulsive connections (Subsection 5.2). Section 6 considers the arbitrary topology case. More precisely, we show convergence to a stable equilibrium almost surely using additional assumptions on the choice of system's parameters. Finally, we apply our stability results to several examples in Section 7, and conclude in Section 8.

## 2 PLASTIC PHASE-COUPLED OSCILLATORS

In this section we first formally describe the model, discuss the meaning of its parameters and define the assumptions that we will use. We then briefly list related works and summarize our results.

### 2.1 Model Description

We study a network of phase-coupled oscillators with plastic coupling strengths whose dynamics are governed by the following two equations:

$$\dot{\phi}_i = \omega_i + \sum_{j \in N_i} K_{ij} \cdot f_{ij}(\phi_j - \phi_i), \ i \in V \tag{1a}$$

$$\dot{K}_{ij} = s_{ij}\big(\alpha_{ij} \cdot (F_{ij}(\phi_j - \phi_i) + q_{ij}) - K_{ij}\big), \ ij \in E, \tag{1b}$$

where $E$ is the set of edges and $V$ the set of vertices. Equation (1a) defines behavior of an oscillator, and equation (1b) determines dynamics of the coupling strength. Here $\phi_i$ is a phase of oscillator $i$ defined on a unit circle $\mathbb{S}^1$ so that all $n$ phase variables are defined on a $n$-dimensional torus $\mathbb{T}^n$; $\omega_i$ is its intrinsic frequency; $N_i$ is a set of oscillators connected to oscillator $i$, i.e. the set of its neighbors; $K_{ij}$ and $f_{ij}$ are a coupling strength and a coupling function, respectively, between connected oscillators $i$ and $j$.

The positive constants $s_{ij}$ in equation (1b) define the rate of change of the coupling strengths, and $F_{ij}(x) \triangleq -\int_0^x f_{ij}(t)\, dt + C_{ij}$ with a choice of integration constant $C_{ij}$ that makes $\int_0^\pi F_{ij}(t)\, dt = 0$. The parameters $q_{ij} \in (-\infty, +\infty)$ and $\alpha_{ij} > 0$ determine the interval of values that the coupling strength $K_{ij}$ can take. More precisely, $K_{ij} \in [\alpha_{ij} \cdot (F_{ij}^{min} + q_{ij}), \ \alpha_{ij} \cdot (F_{ij}^{max} + q_{ij})]$ where $F_{ij}^{min} \leq 0$ and $F_{ij}^{max} \geq 0$ are the minimum and maximum values of the function $F_{ij}$, respectively. Two values of $q_{ij}$ are of a special interest: $q_{ij}^+ \triangleq -F_{ij}^{min} \geq 0$ and $q_{ij}^- \triangleq -F_{ij}^{max} \leq 0$. If $q_{ij} \geq q_{ij}^+$, then the coupling between oscillators $i$ and $j$ is *positive* ($K_{ij} \geq 0$), while the coupling is *negative* ($K_{ij} \leq 0$) when $q_{ij} \leq q_{ij}^-$. In particular, when $q_{ij} = q_{ij}^+ = -F_{ij}^{min}$ (resp. $q_{ij} = q_{ij}^- = -F_{ij}^{max}$), then the coupling strength $K_{ij}$ takes values from the interval $[0, \ \alpha_{ij} \cdot (F_{ij}^{max} - F_{ij}^{min})]$ (resp. $[\alpha_{ij} \cdot (F_{ij}^{min} - F_{ij}^{max}), \ 0]$).

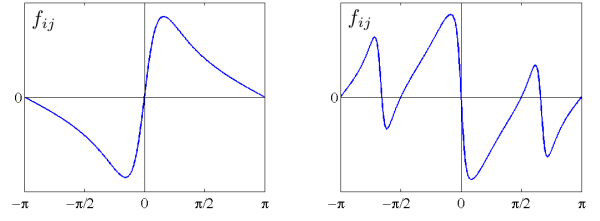

Fig. 1: Examples of functions $f_{ij}$ satisfying Assumption 1.

The topology of system (1) is defined by an undirected connected graph $G = (V, E)$. Each vertex $i \in V$ corresponds to the oscillator $\phi_i$, and each edge $ij \in E$ corresponds to the coupling strength $K_{ij}$, so that $|V| = n$, where $n$ is a number of oscillators in a system, and $N_i = \{j \in V | ij \in E\}$. Additionally, if oscillators $i$ and $j$ are *not* connected, then the coupling strength between them is always equal to zero, i.e. $K_{ij} \equiv 0$ if $ij \notin E$. We denote by $m$ the number of edges in a graph so that $|E| = m$. Therefore, the total number of variables and equations in system (1) is $n+m$. It is assumed that coupling is symmetric, and $K_{ij}$ and $K_{ji}$ are the same variable.

System (1) is fairly general and includes other models studied in the literature as special cases. In particular, when $q_{ij} = 0 \ \forall \ i,j$ and $f_{ij}(\phi_j - \phi_i) = \sin(\phi_j - \phi_i) \ \forall i,j$, then $F_{ij}(\phi_j - \phi_i) = \cos(\phi_j - \phi_i)$ and system (1) becomes the Kuramoto model with varying coupling strengths also known as generalized Kuramoto model [41]. We do require however that the functions $f_{ij}$ satisfy the following three conditions:

**Assumption 1** *Functions $f_{ij} \ \forall ij \in E$ satisfy:*

1) *Symmetric coupling: $f_{ij} = f_{ji}$;*
2) *Odd: $f_{ij}(x) = -f_{ij}(-x)$;*
3) *$C^1$: $f_{ij}$ is continuously differentiable.*

Examples of two functions $f_{ij}$ satisfying these three conditions are shown in Fig. 1.

In this article we study frequency synchronization of system (1). We say that system (1) achieves frequency synchronization if $\dot{\phi}_1(t) = \cdots = \dot{\phi}_n(t) = \dot{\phi}$ and $\dot{K}_{ij}(t) = 0 \ \forall i, j$ as $t \to \infty$, where $\dot{\phi}$ is a common synchronization frequency. In the case of homogeneous oscillators, all intrinsic frequencies of oscillators are equal, i.e. there exists a constant $\omega$ such that $\omega_1 = \cdots = \omega_n = \omega$. It is easy to see that if such homogeneous oscillators synchronize, then their synchronization frequency is $\omega$. Without loss of generality, we assume that $\omega = 0$, and in the rest of the article consider the following system:

$$\dot{\phi}_i = \sum_{j \in N_i} K_{ij} \cdot f_{ij}(\phi_j - \phi_i), \ i \in V \tag{2a}$$

$$\dot{K}_{ij} = s_{ij}\big(\alpha_{ij} \cdot (F_{ij}(\phi_j - \phi_i) + q_{ij}) - K_{ij}\big), \ ij \in E. \tag{2b}$$

Observe that if $(\phi^*, K^*)$ is an equilibrium of system (2), then $(\phi^* + \delta \mathbf{1}_n, K^*)$, where $\mathbf{1}_n$ is a n-dimensional vector of ones and $\delta \in \mathbb{R}$, is also an equilibrium and belongs to the same limit cycle. We will not differentiate between equilibria belonging to the same orbit and thus consider them to be identical. Therefore, in the rest of the article when we talk about the stability of an equilibrium $(\phi^*, K^*)$, we imply stability of the following set of equilibria: [1]

$$E_{\phi^*} = \{(\phi, K) \colon (\phi, K) = (\phi^* + \delta \mathbf{1}_n, K^*), \delta \in \mathbb{R}\}. \quad (3)$$

Further, two equilibria $(\hat{\phi}, \hat{K})$ and $(\bar{\phi}, \bar{K})$ of the same system of plastic phase-coupled oscillators are called topologically equivalent, if they are characterized by the same phase differences, i.e. if $(\hat{\phi}_i - \hat{\phi}_j) = (\bar{\phi}_i - \bar{\phi}_j) \; \forall ij \in E$, or if all phase differences are opposite in sign, i.e. when $(\hat{\phi}_i - \hat{\phi}_j) = -(\bar{\phi}_i - \bar{\phi}_j) \; \forall ij \in E$.

### 2.2 Related Work and Contributions

The plastic phase-coupled oscillator model (1) was initially introduced in [20] as an extension to the classical Kuramoto model to capture the behavior of neural networks. Because the strength of synapses – connections between neurons – can generally change its value and is believed to play a key role in learning and memory formation in the brain, it is natural to consider plastic coupling strengths between oscillators in the Kuramoto model. A well-known synaptic plasticity mechanism called Hebbian rule [18] states that a synapse between two simultaneously active neurons, i.e. neurons that spike almost at the same time, becomes stronger. When neurons are modeled by phase-coupled oscillators, simultaneously firing neurons can be represented by oscillators whose phases are almost equal. This idea is implied in the model (1) with $F_{ij} = \cos()$, where a connection between two oscillators becomes stronger if it is small enough and if the phases of these oscillators are close to each other.

Previous works have introduced and investigated several modifications to (1). For example, in [42], [1], [2] time delays are considered, and the behavior of the system for different values of delay parameters is experimentally explored. In [29] the coupling strength equation of (1) was replaced by an exponential Spike Timing-Dependent Plasticity (STDP) rule, in which a coupling strength's $K_{ij}$ dependence on a phase difference $(\phi_j - \phi_i)$ is defined via exponential function instead of function $F_{ij}$. In [22] a stochastic model of oscillators is studied, where equations (1) contain additive Gaussian noise terms. Synchronization of model (1) with the complete topology, $\sin()$ coupling and $q_{ij} = 0$ is explored in [17] for both, homogeneous and heterogeneous oscillators. While in our previous work [13] results were obtained for model (1) with $q_{ij} = 0 \; \forall \; ij \in E$, arbitrary choices of parameter $q_{ij}$ are considered in this work.

The contributions of our work with respect to the existing literature are manifold. Firstly, we perform a thorough analytical analysis of the system of plastic phase-coupled oscillators in contrast to the empirical studies [8], [35], [39]. Secondly, we consider a fairly general form of the system: instead of studying plastic coupling based on a trigonometric $\sin()$ [40], [41], [42], [17], we investigate a general class of coupling functions. We further explore the

1. Alternatively, we could consider phase differences $(\phi_i - \phi_j)$ as the variables of system (2) and study stability of a single equilibrium instead of a set (3) of equilibria. This approach will be used in Section 6.

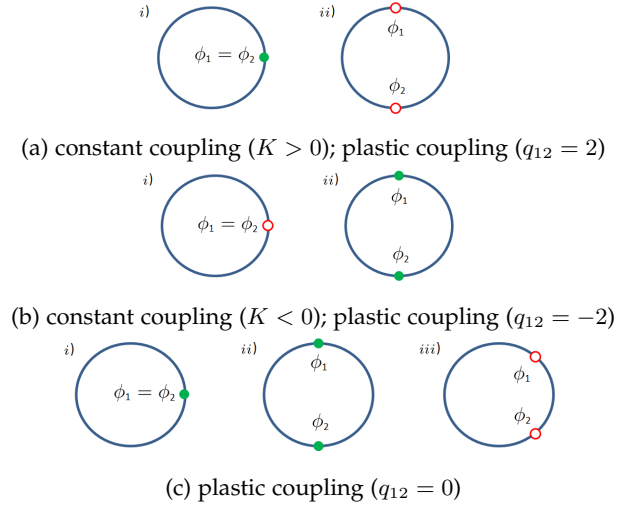

(a) constant coupling ($K > 0$); plastic coupling ($q_{12} = 2$)

(b) constant coupling ($K < 0$); plastic coupling ($q_{12} = -2$)

(c) plastic coupling ($q_{12} = 0$)

Fig. 2: Example with two oscillators. Green (filled) circles correspond to stable equilibria, red (not filled) circles correspond to unstable equilibria.

behavior of model (2) for various values of its parameters and show how they impact the properties of the model such as synchronization and stability. Finally, several interesting and important examples are provided that illustrate specific features of the system and confirm our theoretical results.

## 3 MOTIVATING EXAMPLES

In this section we illustrate the differences between the plastic oscillator model (2) and the constant coupling counter part using several simple examples. The number of oscillators in these examples varies from two to four, and, for illustrative purpose, in this section we assume that $f_{ij}(\phi_j - \phi_i) = \sin(\phi_j - \phi_i)$ for all connected oscillators $i$ and $j$. We demonstrate with these examples that stability of equilibria may change if the coupling strength becomes plastic. Additionally, new equilibria may arise in this case including scenarios with infinitely many equilibria as shown in Subsection 3.3. Furthermore, the set of equilibria points and their stability may depend on the value of the parameter $q_{ij}$. For each example we consider three types of coupling strengths: constant equal positive, constant equal negative, and varying coupling strengths. For the varying or plastic coupling strength we additionally explore the cases of *strictly positive* ($q_{ij} > q_{ij}^+ = 1$), *strictly negative* ($q_{ij} < q_{ij}^- = -1$) and *symmetric hybrid* ($q_{ij} = 0$) connections. All examples presented here will be used again for illustrative purposes in Section 7, where we apply our theoretical results to each example to explore its dynamics and stability of its equilibria.

### 3.1 Two Oscillators

A system of two connected homogeneous oscillators with a constant coupling strength $K$ and a trigonometric $\sin()$ coupling function is described by the following equations:

$$\dot{\phi}_i = K \cdot \sin(\phi_j - \phi_i), \quad (4)$$

where $i = 1$, $j = 2$ or $i = 2$, $j = 1$. This system has two topologically distinct equilibria: one is in-phase and stable: $\phi_1 = \phi_2$ (Fig. 2a-i), and the other one is anti-phase and unstable: $\phi_2 = \phi_1 + \pi$ (Fig. 2a-ii). When the coupling strength $K$ is constant and negative, then the set of equilibria of system (4) remains unchanged, but the in-phase equilibrium now becomes unstable (Fig. 2b-i), whereas the anti-phase equilibrium becomes stable (Fig. 2b-ii).

System (2) of two oscillators with a plastic coupling strength and $q_{12} = 0$ contains two sets of equilibria that are characterized by the following conditions:

1). $\sin(\phi_2 - \phi_1) = 0$. If $\phi_1 = \phi_2$, then $K = \alpha$, and when $\phi_2 = \phi_1 + \pi$, then $K = -\alpha$.

2). $K = 0$, then $\cos(\phi_2 - \phi_1) = 0$, i.e. $\phi_2 = \phi_1 + \pi/2$. The Jacobian for the system of two oscillators takes form:

$$J = \begin{bmatrix} -K \cdot \cos() & K \cdot \cos() & \sin() \\ K \cdot \cos() & -K \cdot \cos() & -\sin() \\ \alpha \cdot s \cdot \sin() & -\alpha \cdot s \cdot \sin() & -s \end{bmatrix},$$

where $\sin() = \sin(\phi_2 - \phi_1)$ and $\cos() = \cos(\phi_2 - \phi_1)$ for brevity. It can be easily verified that equilibria from condition 1) above are stable (Fig. 2c-i, 2c-ii), whereas the equilibria from condition 2) (when $K = 0$) are unstable (Fig. 2c-iii).

If the coupling strength is plastic and $q_{12} = 2$ (*positive* coupling), then the set of equilibria and their stability are the same as in the case of a constant positive $K$ (Fig. 2a). Similarly, if $q_{12} = -2$ (*negative* coupling), the equilibria and stability coincide with the ones corresponding to the case of a constant negative coupling $K$ (Fig. 2b). This property is a priori not necessarily true given the fact that system (2) has a larger state space that can in principle change the stability of an equilibrium.

Therefore, two observations can be made: all equilibria of system (4) with constant coupling strength are also equilibria of system (2) with plastic coupling strength for each of three values of $q_{12}$. The second observation is that when the coupling strength is non-constant and $q_{12} = 0$, a new set of equilibria emerges. This set, however, contains only unstable equilibria in the case of two oscillators.

### 3.2 Three Oscillators

In this subsection we consider an example of three connected homogeneous oscillators. We assume that underlying topology is a complete graph, which means that each oscillator is connected to two others. When the coupling function is $\sin()$ and coupling strength $K$ is constant, the behavior of the oscillators is defined by the following set of equations:

$$\begin{aligned} \dot\phi_1 &= K \cdot \sin(\phi_2 - \phi_1) + K \cdot \sin(\phi_3 - \phi_1), \\ \dot\phi_2 &= K \cdot \sin(\phi_1 - \phi_2) + K \cdot \sin(\phi_3 - \phi_2), \\ \dot\phi_3 &= K \cdot \sin(\phi_1 - \phi_3) + K \cdot \sin(\phi_2 - \phi_3). \end{aligned} \quad (5)$$

When $K > 0$, system (5) has 3 topologically distinct equilibria (Fig. 3): one is in-phase when $\phi_1 = \phi_2 = \phi_3$ and stable (Fig. 3a-i), another one is $\phi_1 = \phi_2$, $\phi_3 = \phi_1 + \pi$ and unstable (Fig. 3a-ii), and the last one is defined as $\phi_2 = \phi_1 + 2\pi/3$, $\phi_3 = \phi_1 - 2\pi/3$ and is also unstable (Fig. 3a-iii).

When $K < 0$, the set of equilibria of system (5) remains the same. Stability properties of the equilibria, however,

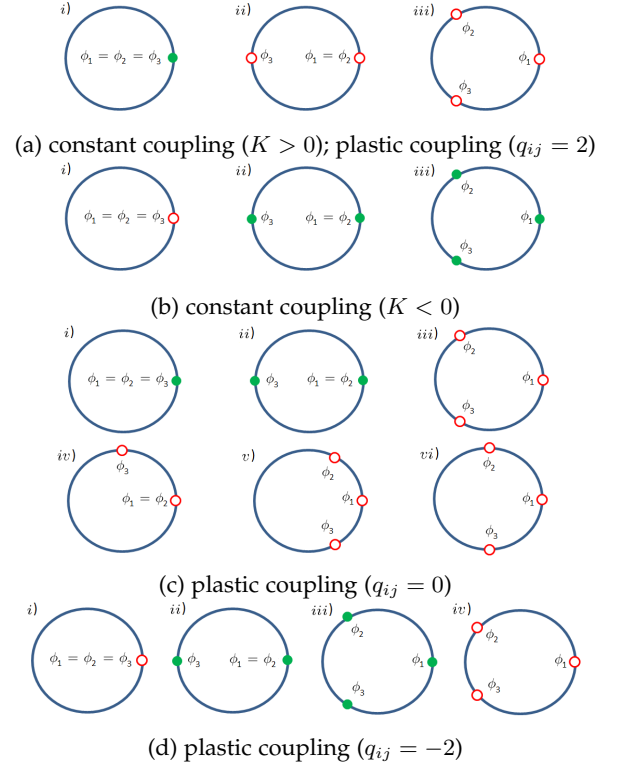

(a) constant coupling ($K > 0$); plastic coupling ($q_{ij} = 2$)

(b) constant coupling ($K < 0$)

(c) plastic coupling ($q_{ij} = 0$)

(d) plastic coupling ($q_{ij} = -2$)

Fig. 3: Example with three oscillators. Green (filled) circles correspond to stable equilibria, red (not filled) circles correspond to unstable equilibria.

change, as in the example with two oscillators. In particular, the first equilibrium becomes unstable (Fig. 3b-i), and the last two become stable (Fig. 3b-ii, 3b-iii).

We now consider the case of plastic coupling strengths while assuming that $\alpha_{12} = \alpha_{23} = \alpha_{13} = \alpha > 0$. When $q_{ij} = 0$ $\forall ij \in E$, a detailed description of equilibria was provided in [13] and is omitted here. The equilibria and their stability when $q_{ij} = 0$ $\forall ij \in E$ are illustrated in Fig. 3c.

When $q_{ij} = 2$ for each edge $ij \in E$, equilibria and their stability coincide with the case of constant positive coupling strength (Fig. 3a), but when $q_{ij} = -2$ ($\forall ij \in E$), a new equilibrium emerges (Fig. 3d-iv). At this unstable equilibrium, $(\phi_2 - \phi_1) = (\phi_1 - \phi_3) \approx 0.785\pi$, $K_{12} = K_{13} \approx -2.7808\alpha$, and $K_{23} \approx -1.7808\alpha$.

From the considered examples of two and three oscillators several observations can be made. First, each equilibrium of a system with constant coupling strengths was also an equilibrium of the corresponding system with varying coupling strengths. This is not true, however, for all values of parameters $\alpha$ and $q$. Second, stability of these common equilibria can differ for systems with constant and non-constant coupling. Third, system (2) can possess additional equilibria, and moreover, the set of equilibria and their stability may depend on the parameter $q$.

### 3.3 Four Oscillators

We consider here the case of four oscillators connected by a complete graph. Instead of describing all equilibria of this system, we will show that system (2) with four homogeneous oscillators, $\sin()$ coupling, equal $\alpha_{ij} = \alpha > 0$ $\forall ij \in E$, and $q_{ij} = 0$ $\forall ij \in E$, has infinitely many topologically distinct equilibria.

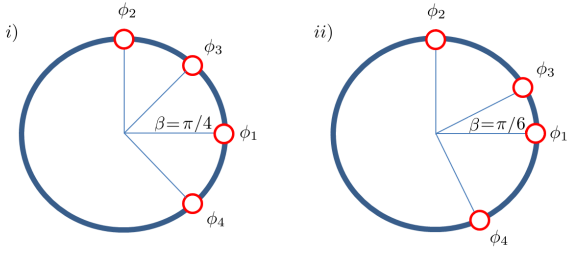

Fig. 4: Equilibria corresponding to $\beta = \pi/4$ (left) and $\beta = \pi/6$ (right) for the example with four oscillators.

These equilibria can be defined by means of a parameter $\beta$. Then, for each value of $\beta \in (0, \pi/2)$, phases: $\phi_2 = \phi_1 + \pi/2$, $\phi_3 = \phi_1 + \beta$, $\phi_4 = \phi_3 - \pi/2$, and coupling strengths $K_{12} = K_{34} = 0$, $K_{13} = \cos(\beta) \cdot \alpha$, $K_{14} = K_{23} = \cos(\pi/2 - \beta) \cdot \alpha = \sin(\beta) \cdot \alpha$, $K_{24} = \cos(\pi - \beta) \cdot \alpha = -\cos(\beta) \cdot \alpha$ define an equilibrium. Phases corresponding to this equilibrium with values of parameter $\beta = \pi/4$ and $\beta = \pi/6$ are shown in Fig. 4. Notice, that in all such equilibria two coupling strengths are equal to zero, and edges corresponding to non-zero coupling strengths form a graph with a ring topology.

The infinite set of equilibria defined for the case of four oscillators can be generalized for all systems with even $n > 2$ number of oscillators: $\phi_1 \ldots \phi_{\frac{n}{2}-1}$ have the same phase $\phi$, oscillators $\phi_{\frac{n}{2}} \ldots \phi_{n-2}$ have phase $\phi - \frac{\pi}{2}$, and two other oscillators have phases $\phi_{n-1} = \phi - \beta$ and $\phi_n = \phi - \beta - \frac{\pi}{2}$.

## 4 SYNCHRONIZATION AND STABILITY ANALYSIS

This section contains general results obtained for system (2) with arbitrary values of parameters $q_{ij}$ and arbitrary underlying topology. We first show in Theorem 1 by providing a Lyapunov function that system (2) of homogeneous oscillators is gradient, and thus always converges to a set of equilibria, i.e. achieves frequency synchronization (Subsection 4.1). After that we formulate sufficient instability and stability conditions for equilibria of system (2) with arbitrary underlying topology and arbitrary $q_{ij}$ in Theorems 2 and 3, respectively (Subsection 4.2). We then apply the derived results to explore the stability of the in-phase and anti-phase equilibria in Subsection 4.3.

### 4.1 Frequency Synchronization

In Theorem 1 we prove that system (2) of homogeneous oscillators is a gradient system and always achieves frequency synchronization. A similar result was obtained in [40], where a potential function was found for system (2) and frequency synchronization was shown, but only for the case of a complete graph topology and for $q_{ij} = 0, \forall ij \in E$. Therefore, Theorem 1 generalizes the result from [40] for the case of an arbitrary topology and arbitrary values of parameters $q_{ij}$.

**Theorem 1 (Frequency synchronization)** *System* (2) *is a gradient system and achieves frequency synchronization for all initial values of phases and coupling strengths.*

*Proof.* Notice that all phase variables $\phi_i$ are defined on a n-dimensional torus $\mathbb{T}^n$ which is compact. The result of the theorem holds if the coupling strengths $K_{ij}$ are defined on the whole $\mathbb{R}^m$. Indeed, we can provide a potential function $V$:

$$V = - \sum_{ij \in E, i>j} K_{ij} F_{ij}(\phi_i - \phi_j) - \sum_{ij \in E, i>j} q_{ij} K_{ij} + \frac{1}{2} \sum_{ij \in E, i>j} \frac{K_{ij}^2}{\alpha_{ij}}. \quad (6)$$

This function is well-defined, radially unbounded, bounded below, and it is easy to verify that the derivative of $V$ with respect to time is:

$$\dot{V} = (\nabla V)^T * \begin{bmatrix} \dot{\phi}_i \\ . \\ . \\ \dot{\phi}_n \\ \dot{K}_{21} \\ . \\ . \\ \dot{K}_{n,n-1} \end{bmatrix} = -\sum_{i=1}^{n} (\dot{\phi}_i)^2 - \sum_{ij \in E, i>j} \frac{(\dot{K}_{ij})^2}{\alpha_{ij} s_{ij}}.$$

We can see that $\dot{V}$ is always non-positive and is equal to zero if and only if $\dot{\phi}_i = 0$ and $\dot{K}_{ij} = 0$ for all $i$ and $j$. Thus, by LaSalle's Invariance Principle [27], the trajectories of (2) always converge to a set of equilibria. In other words, for all initial conditions frequency synchronization occurs. □

**Remark** *Notice that Theorem 1 does not imply pointwise convergence to a single equilibrium. It is also not guaranteed that equilibria of system* (2) *are isolated.*

### 4.2 Stability and Instability Conditions

The main results of this subsection are Theorem 2 which is a sufficient instability condition, and Theorem 3 that defines a sufficient condition for stability of an equilibrium of system (2). These results are based on Lyapunov's indirect method [24], that states:

1) If $\text{Re}\lambda_i < 0$ for all eigenvalues of the Jacobian matrix $J$, then equilibrium is asymptotically stable.
2) If $\text{Re}\lambda_i > 0$ for at least one eigenvalue of the Jacobian matrix $J$, then equilibrium is unstable.

Let $B \in \mathbb{R}^{n \times m}$ denote an oriented incidence matrix of a graph that defines underlying topology of system (2). Then element $(i, e)$ of this matrix is

$$B(i, e) = \begin{cases} 1 & \text{if } i \text{ is the head of } e, \\ -1 & \text{if } i \text{ is the tail of } e, \\ 0 & \text{otherwise}, \end{cases} \quad (7)$$

where $e$ is an edge of graph $G$. Although the definition of matrix $B$ implies that $G$ is oriented, all properties of this matrix used in this article do not depend on a particular orientation. Therefore, we assume that for a given undirected graph $G$, an arbitrary orientation of its edges is chosen, i.e. for every undirected edge $e = ij$ one of the nodes $i, j$ is designated as the head of $e$, and another one corresponds to the tail of $e$.

Let $\alpha \in \mathbb{R}^m$, $s \in \mathbb{R}^m$, $K \in \mathbb{R}^m$, $f' \in \mathbb{R}^m$ and $f \in \mathbb{R}^m$ denote vectors whose components are $\alpha_{ij}$, $s_{ij}$, $K_{ij}$, $f'_{ij}(\phi_j - \phi_i)$ and $f_{ij}(\phi_j - \phi_i)$, respectively, for each $i, j$ such that $ij \in E$. We will use symbol $*$ to denote the componentwise product of vectors. Jacobian of system (2) does not depend on the values of parameters $q_{ij}$ and can be written in a following way:

$$J = \begin{bmatrix} B & 0 \\ 0 & I \end{bmatrix} \begin{bmatrix} -\text{diag}(K*f') & -\text{diag}(f) \\ -\text{diag}(\alpha*s*f) & -\text{diag}(s) \end{bmatrix} \begin{bmatrix} B^T & 0 \\ 0 & I \end{bmatrix}.$$

The first matrix in the product is of size $(n+m) \times (m+m)$, the second matrix is of size $(m+m) \times (m+m)$ and the last matrix in the product has dimensions $(m+m) \times (n+m)$. Notice that Jacobian $J$ has a trivial eigenvector $[\mathbf{1}_n\ \mathbf{0}_m]^T$, where $\mathbf{1}_n \in \mathbb{R}^n$ is a vector of ones and $\mathbf{0}_m \in \mathbb{R}^m$ is a vector of zeros with $n$ and $m$ components, respectively. This eigenvector emerges due to rotational invariance of system (2) and corresponds to a zero eigenvalue. Since trajectories of system (2) are orthogonal to the direction of an orbit, we still can apply Lyapunov's indirect method to explore stability of the set (3). If all remaining eigenvalues of the Jacobian have negative real part, then equilibrium is stable; if there exists an eigenvalue with a positive real part, then equilibrium is unstable.

The component of vector $f$ that corresponds to the edge $e = ij$ is equal to $f_{ij}(\phi_i - \phi_j)$ if edge $e = ij$ is oriented from a tail $j$ to a head $i$, and thus $B(i, e) = 1$, $B(j, e) = -1$. Similarly, if edge $e = ij$ is oriented from a tail $i$ to a head $j$, then $B(i, e) = -1$, $B(j, e) = 1$ and component of $f$ associated with edge $ij$ is equal to $f_{ij}(\phi_j - \phi_i)$.

Each partition $P$ of the graph's vertices into two sets $V^-$ and $V^+$ such that $V^- \cap V^+ = \varnothing$ and $V^- \cup V^+ = V$, defines a cut $C(P) \triangleq \{ij \in E | i \in V^-, j \in V^+\}$. With each cut $C(P)$ we associate a cut vector $c_P \in \mathbb{R}^m$ which is defined as follows:

$$c_P(e) = \begin{cases} 1 & \text{if } e \text{ goes from } V^- \text{ to } V^+, \\ -1 & \text{if } e \text{ goes from } V^+ \text{ to } V^-, \\ 0 & \text{if } e \notin C(P). \end{cases} \quad (8)$$

We can now formulate the following instabiliy condition that is similar to Theorem 2 of [30]:

**Theorem 2 (Sufficient instability condition)** *If there exists a cut $C(P)$ such that at equilibrium $(\phi^*, K^*)$ of system (2):*

$$\sum_{ij \in C(P)} (K_{ij} f'_{ij} - \alpha_{ij} f_{ij}^2) < 0, \quad (9)$$

*where $K_{ij} = K^*_{ij}$, $f'_{ij} = f'_{ij}(\phi^*_j - \phi^*_i)$ and $f_{ij} = f_{ij}(\phi^*_j - \phi^*_i)$, then $(\phi^*, K^*)$ is an unstable equilibrium.*

*Proof.* We first show that the Jacobian of system (2) can be decomposed into a product of matrices $D$ and $A$:

$$J = DA, \quad (10)$$

where $D$ is a positively-definite diagonal matrix, and $A$ is a symmetric matrix. We then demonstrate that stability of equilibria of system (2) does not depend on matrix $D$, because matrices $J$ and $A$ have the same number of positive, negative and zero eigenvalues. Next, for matrix $A$ we provide a vector $\vec{X}$ such that $\vec{X}^T A \vec{X} > 0$, which guarantees that the symmetric matrix $A$ has a positive eigenvalue and so does the Jacobian matrix $J$. This in turn means that an equilibrium is unstable due to Lyapunov's indirect method.

Decomposition (10) is possible because system (2) is a gradient system. Note that the Hessian matrix $H(V)$ of the potential function $V$ is symmetric. Let diagonal $(n+m) \times (n+m)$ matrix $D$ be defined as

$$D = \begin{bmatrix} I & 0 \\ 0 & \text{diag}(\alpha*s) \end{bmatrix}, \quad (11)$$

then, since equations (2) can be written as follows:

$$\begin{bmatrix} \dot{\phi} \\ \dot{K} \end{bmatrix} = -D \cdot \nabla V, \quad (12)$$

decomposition (10) exists with $A = -H(V)$.

We now show that matrices $J$ and $A = -H(V)$ have the same numbers of positive, negative and zero eigenvalues. Observe that if matrix $D^{\frac{1}{2}}$ is a square root of matrix $D$, then matrices $DA$ and $D^{\frac{1}{2}} A D^{\frac{1}{2}}$ have the same eigenvalues, because matrix $D$ is positive-definite. This also implies that Jacobian of system (2) with homogeneous oscillators has only real eigenvalues. Next, since $A$ is a symmetric matrix with real entries, it can be diagonalized by an orthogonal matrix, i.e. there exists a real orthogonal matrix $Q$ such that $A = QGQ^T$, where $G$ is a diagonal matrix. Further, notice that

$$D^{\frac{1}{2}} A D^{\frac{1}{2}} = D^{\frac{1}{2}} Q G Q^T D^{\frac{1}{2}} = LGL^T, \quad (13)$$

where matrix $L$ is defined as $L = D^{\frac{1}{2}} Q$ and is invertible. Therefore,

$$Q^T A Q = L^{-1}(D^{\frac{1}{2}} A D^{\frac{1}{2}})(L^{-1})^T = G. \quad (14)$$

By Sylvester's law of inertia [7], numbers of positive, negative and zero eigenvalues of matrices $A$, $D^{\frac{1}{2}} A D^{\frac{1}{2}}$ and $G$ are equal. Thus, since $J = DA$ and $D^{\frac{1}{2}} A D^{\frac{1}{2}}$ have equal eigenvalues, then the numbers of positive, negative and zero eigenvalues of matrices $J$ and $A$ are the same.

We now consider the symmetric matrix $A$ and show that when condition (9) is satisfied, matrix $A$ has a positive eigenvalue. We define a symmetric $(2m) \times (2m)$ matrix $M$ to be:

$$M = \begin{bmatrix} \text{diag}(K*f') & \text{diag}(f) \\ \text{diag}(f) & \text{diag}(1/\alpha) \end{bmatrix}, \quad (15)$$

where $1/\alpha$ is a vector with components $1/\alpha_{ij}$, then

$$A = -\begin{bmatrix} B & 0 \\ 0 & I \end{bmatrix} M \begin{bmatrix} B^T & 0 \\ 0 & I \end{bmatrix}. \quad (16)$$

We denote:

$$\hat{B} = \begin{bmatrix} B & 0 \\ 0 & I \end{bmatrix}, \quad (17)$$

and then

$$A = -\hat{B} M \hat{B}^T. \quad (18)$$

Now we will assume that there exists a cut $C(P)$ that satisfies condition (9). We define a vector $\vec{Y} \in R^{2m}$ to be:

$$\vec{Y} = \begin{bmatrix} c_P \\ -c_P * f * \alpha \end{bmatrix}, \quad (19)$$

where $c_P$ is a cut vector associated with the cut $C(P)$, and multiplication in $-c_P * f * \alpha$ is componentwise.

It can be verified that the sum from (9) is equal to $\vec{Y}^T M \vec{Y}$. Indeed, if an edge $k$ ($1 \leq k \leq m$) belongs to the cut $C(P)$, then $\vec{Y}_k = \pm 1$ and $\vec{Y}_{k+m} = \mp f_k \alpha_k$. The summand number $k$ in $\vec{Y}^T M \vec{Y}$ is equal to:

$$Y_k^2 K_k f_k' + 2 Y_k Y_{m+k} f_k + \frac{Y_{m+k}^2}{\alpha_k} \qquad (20)$$
$$= K_k f_k' - 2 f_k^2 \alpha_k + f_k^2 \alpha_k = K_k f_k' - f_k^2 \alpha_k,$$

which is also the $k^{th}$ summand of the sum (9).

The cut space of the graph $G$ is defined as a space spanned by all cut vectors $c_P$. It is known (see for example [4]) that the range of $B^T$ is the cut space of $G$. Therefore, for the cut vector $c_P$ there exists a vector $\vec{x}_1 \in R^n$ such that $c_P = B^T \vec{x}_1$. Therefore,

$$\vec{Y} = \begin{bmatrix} B^T \vec{x}_1 \\ -c_P * f * \alpha \end{bmatrix} = \begin{bmatrix} B^T & 0 \\ 0 & I \end{bmatrix} \cdot \begin{bmatrix} \vec{x}_1 \\ -c_P * f * \alpha \end{bmatrix} = \hat{B}^T \vec{X}, \qquad (21)$$

where $\vec{X} = \begin{bmatrix} \vec{x}_1 \\ -c_P * f * \alpha \end{bmatrix} \in R^{n+m}$. Finally,

$$0 > \vec{Y}^T M \vec{Y} = \vec{X}^T \hat{B} M \hat{B}^T \vec{X} = -\vec{X}^T A \vec{X}, \qquad (22)$$

which means that there is a vector $\vec{X}$ such that $\vec{X}^T A \vec{X} > 0$ and thus symmetric matrix $A$ has a positive eigenvalue which implies that Jacobian $J$ has also a positive eigenvalue. Therefore, equilibrium $(\phi^*, K^*)$ is unstable. $\square$

We now formulate a sufficient condition for an equilibrium of system (2) to be stable.

**Theorem 3 (Sufficient stability condition)** *If at equilibrium $(\phi^*, K^*)$ of system (2), for each $ij \in E$:*

$$K_{ij} f_{ij}' - \alpha_{ij} f_{ij}^2 > 0, \qquad (23)$$

*where $K_{ij} = K_{ij}^*$, $f_{ij}' = f_{ij}'(\phi_j^* - \phi_i^*)$ and $f_{ij} = f_{ij}(\phi_j^* - \phi_i^*)$, then equilibrium $(\phi^*, K^*)$ is asymptotically stable.*

*Proof.* All eigenvalues of the Jacobian of system (2) are real. To apply Lyapunov's indirect method, we need to show that at equilibrium $(\phi^*, K^*)$ Jacobian has only negative eigenvalues. However, it has always at least one zero eigenvalue that corresponds to the rotational invariance of the system: if all phases $\phi_i$ ($i = 1, \ldots, n$) are simultaneously shifted by the same value, the system does not change. The eigenvector associated with this zero eigenvalue is a vector $[\mathbf{1}_n \ \mathbf{0}_m]^T$. As previously mentioned, in this article we do not distinguish equilibria that belong to the same set (3), and thus study stability of the whole set $E_{\phi^*}$. To show stability of $E_{\phi^*}$ using an indirect Lyapunov's method, we need to show that all remaining eigenvalues of the Jacobian are strictly negative.

In the proof of Theorem 2 it was shown that the Jacobian matrix $J$ and symmetric matrix $A$ have the same numbers of negative, positive and zero eigenvalues. This means that matrix $A$ also possesses a zero eigenvalue corresponding to the rotational invariance. Moreover, it is easy to see that vector $[\mathbf{1}_n \ \mathbf{0}_m]^T$ is also an eigenvector of matrix $A$ associated with a zero eigenvalue. Therefore, to prove that equilibrium $(\phi^*, K^*)$ is stable, it is sufficient to demonstrate that all eigenvalues of matrix $A$ are negative (except for one zero eigenvalue corresponding to the rotational invariance), or that $\vec{X}^T A \vec{X} < 0$ for all non-zero vectors $\vec{X} \in \mathbb{R}^{n+m}$, $\vec{X} \notin \mathrm{span}([\mathbf{1}_n \ \mathbf{0}_m]^T)$, since $A$ is symmetric.

Notice that because $A = -\hat{B} M \hat{B}^T$, the matrix $A$ will have only negative eigenvalues (except one) if $\vec{Y}^T M \vec{Y} > 0$ for all non-zero vectors $\vec{Y} \in \mathbb{R}^{2m}$. Indeed, if $\hat{B}^T \vec{X} \neq \mathbf{0}_{n+m}$, then

$$\vec{X}^T A \vec{X} = -\vec{X}^T \hat{B} M \hat{B}^T \vec{X} = -\vec{Y}^T M \vec{Y} < 0, \qquad (24)$$

where the vector $\vec{Y} \triangleq \hat{B}^T \vec{X}$.

Additionally, if $\vec{X} = [\vec{x}_1 \ \vec{x}_2]^T$, where $\vec{x}_1$ are the first $n$ components of $\vec{X}$, and $\vec{x}_2$ are the last $m$ components of $\vec{X}$, then $\hat{B}^T \vec{X} = \mathbf{0}_{n+m}$ only if $B^T \vec{x}_1 = \mathbf{0}_n$ and $\vec{x}_2 = \mathbf{0}_m$. And since $\ker(B^T) = \mathrm{span}(\mathbf{1}_n)$ for a connected $G$ (see for example [4]), then $\hat{B}^T \vec{X} \neq \mathbf{0}_{n+m}$ if $\vec{X} \notin \mathrm{span}([\mathbf{1}_n \ \mathbf{0}_m]^T)$.

Therefore, it is now enough to show that condition (23) is sufficient for matrix $M$ to be positive definite. Let $\vec{Y} \in \mathbb{R}^{2m}$ be an arbitrary vector, then $\vec{Y}^T M \vec{Y}$ is a sum of $m$ terms, where the $k^{th}$ term is equal to

$$Y_k^2 K_k f_k' + 2 Y_k Y_{m+k} f_k + \frac{Y_{m+k}^2}{\alpha_k}. \qquad (25)$$

We now consider this term as a quadratic function of $Y_k$. This equation is the equation of a parabola whose branches are directed upwards because $K_k f_k' > 0$ due to (23). Then, the minimum value of (25) is achieved at the vertex of the parabola and is equal to:

$$\left(\frac{Y_{m+k} f_k}{K_k f_k'}\right)^2 K_k f_k' - 2\left(\frac{Y_{m+k} f_k}{K_k f_k'}\right) Y_{m+k} f_k + \frac{Y_{m+k}^2}{\alpha_k}$$
$$= -\frac{Y_{m+k}^2 f_k^2}{K_k f_k'} + \frac{Y_{m+k}^2}{\alpha_k} = Y_{m+k}^2 \left(-\frac{f_k^2}{K_k f_k'} + \frac{1}{\alpha_k}\right). \qquad (26)$$

The last expression is positive if $Y_{m+k} \neq 0$ and if condition (23) is satisfied.

Suppose that $Y_{m+k} = 0$, then (25) becomes equal to $Y_k^2 K_k f_k' \geq 0$, and is equal to zero only if $Y_k = 0$. Since $\vec{Y}$ is a non-zero vector, there exists at least one component $k$ of vector $\vec{Y}$ such that the sum (25) is strictly positive, and for all other components these sums are non-negative. Therefore, for all vectors $\vec{Y} \in \mathbb{R}^{2m}$: $\vec{Y}^T M \vec{Y} > 0$, and $\vec{X}^T A \vec{X} < 0$ for all vectors $\vec{X} \in \mathbb{R}^{n+m}$ such that $\vec{X} \notin \mathrm{span}([\mathbf{1}_n \ \mathbf{0}_m]^T)$. Thus, all eigenvalues of $A$ except one are negative, so are eigenvalues of $J$, and therefore th equilibrium $(\phi^*, K^*)$ is asymptotically stable. $\square$

**Remark** *Condition* (23) *is equivalent to the following condition:*

$$(F_{ij} + q_{ij}) \cdot f_{ij}' - f_{ij}^2 > 0, \qquad (27)$$

*where $F_{ij} = F_{ij}(\phi_j^* - \phi_i^*)$, since $K_{ij} = \alpha_{ij}(F_{ij} + q_{ij})$ at an equilibrium, and $\alpha_{ij} > 0$.*

We have proved frequency synchronization of system (2) and found sufficient stability and instability conditions of its equilibria for a fairly general class of functions $f_{ij}$. In the next subsection we will apply these conditions for a more specific class of these functions to investigate stability of in-phase and anti-phase equilibria.

### 4.3 Stability of In-phase and Anti-phase Equilibria

In this subsection we investigate the stability properties of two special types of equilibria of system (2): in-phase and anti-phase equilibria. Equilibrium $(\phi^*, K^*)$ is called in-phase, if $\phi_1^* = \phi_2^* = \cdots = \phi_n^*$, while for an anti-phase equilibrium the absolute value of the phase difference between any two oscillators is either zero or $\pi$: $|\phi_i^* - \phi_j^*| = 0$ or $|\phi_i^* - \phi_j^*| = \pi$ for any $i, j$. To exclude the in-phase equilibrium from the set of anti-phase equilibria, we additionally require that at any anti-phase equilibrium at least for one pair of oscillators $i$ and $j$, their phase difference is equal to $\pi$. Such in-phase and anti-phase states are indeed equilibria of system (2) because $f_{ij}(0) = f_{ij}(\pi) = 0$ for any $i$ and $j$ due to the Assumption 1. Notice that the in-phase equilibrium is unique (up to rotational symmetry), and there are $2^{n-1} - 1$ topologically distinct anti-phase equilibria.

In the rest of the article we concentrate on a more special class of functions $f_{ij}()$. In particular, these functions must fulfill the following conditions.

**Assumption 2** *The functions $f_{ij}$ $\forall ij \in E$ satisfy:*

1) *Assumption 1;*
2) $f'_{ij}(0) > 0$, $f'_{ij}(\pi) < 0$;
3) $f_{ij}(x) > 0$, $\forall x \in (0, \pi)$.

Example of a function that meets all conditions of Assumption 2 is shown on the left side of Fig. 1. Notice, that for instance, function $f_{ij}() = \sin()$ belongs to this type of functions. If $f_{ij}()$ satisfies Assumptions 2, then its corresponding function $F_{ij}$ is strictly decreasing on the interval $[0, \pi]$, and $F_{ij}(0) = F_{ij}^{max} > 0$, $F_{ij}(\pi) = F_{ij}^{min} < 0$. Therefore, there exists a single point $x \in (0, \pi)$ such that $F_{ij}(x) = 0$. This property is crucial for showing isolation of equilibria of system (2) with a tree topology in Corollary 9. If the function $f_{ij}$ satisfies Assumption 2, then *positive* coupling ($q_{ij} \geq q_{ij}^+$) between oscillators $i$ and $j$ is also attractive, and *negative* coupling ($q_{ij} \leq q_{ij}^-$) is repulsive. Thus, in the rest of the article we will use these concepts interchangeably.

For any set of phase values $\phi_1, \ldots, \phi_n$ there exists a unique arc of a circle $\mathbb{S}^1$ that contains all phase values and has a minimum possible length. Let $d(\phi)$ denote the length of this arc for phases $\phi_1, \ldots, \phi_n$. Next theorem provides characterization of the in-phase equilibrium.

**Theorem 4 (Stability of the in-phase equilibrium)** *If functions $f_{ij}()$ satisfy Assumption 2, then the in-phase equilibrium $(\phi^*, K^*)$ is:*

- **asymptotically stable**, *if $q_{ij} > q_{ij}^-$ $\forall ij \in E$;*
- **unstable**, *if $q_{ij} < q_{ij}^-$ $\forall ij \in E$.*

*Moreover, if $q_{ij} < q_{ij}^-$ $\forall ij \in E$, or $q_{ij} > q_{ij}^+$ $\forall ij \in E$, then the in-phase equilibrium is the only equilibrium satisfying $d(\phi) < \pi$.*

*Proof.* The first part of this theorem can be shown by a direct application of Theorems 2 and 3. Indeed, for the in-phase equilibrium and functions $f_{ij}$ satisfying Assumption 2: $f_{ij}(0) = 0$ and $f'_{ij}(0) > 0$ for each $ij \in E$. Further, $(F_{ij}(0) + q_{ij}) > 0$ if $q_{ij} > q_{ij}^-$, because $F_{ij}(0) = F_{ij}^{max}$ (function $F_{ij}$ is decreasing on $[0, \pi]$). Then, $(F_{ij}(0) + q_{ij}) \cdot f'_{ij}(0) > 0$ $\forall ij \in E$, condition (27) is satisfied, and the equilibrium is stable by Theorem 3. Similarly, $(F_{ij}(0) + q_{ij}) < 0$ if $q_{ij} < q_{ij}^-$ which means that $\alpha_{ij}(F_{ij}(0) + q_{ij}) \cdot f'_{ij}(0) = K_{ij} f'_{ij}(0) < 0$ $\forall ij \in E$, condition (9) is satisfied, and the equilibrium is unstable according to Theorem 2.

The second part of this theorem is similar to Lemma 3 of [30] and can be proved as follows. Condition $d(\phi) < \pi$ implies that all phases of oscillators belong to the same half-circle. Suppose that $q_{ij} < q_{ij}^-$ $\forall ij \in E$ (strictly repulsive coupling), or $q_{ij} > q_{ij}^+$ $\forall ij \in E$ (strictly attractive coupling) and there exists a non-in-phase equilibrium $(\hat{\phi}, \hat{K})$ with $d(\phi) < \pi$. Let $\phi_{min}$ be the minimal phase value among all oscillators for this equilibrium. Then at least one oscillator $k$ with phase value $\hat{\phi}_k = \phi_{min}$ is connected to an oscillator $l$ with strictly greater phase value $\hat{\phi}_l > \hat{\phi}_k$, because the graph is assumed to be connected. Therefore, $0 < (\hat{\phi}_l - \hat{\phi}_k) < \pi$ and $f_{kl}(\hat{\phi}_l - \hat{\phi}_k) > 0$ due to the Assumption 2. If $q_{ij} < q_{ij}^-$ $\forall ij \in E$, then all coupling strengths are strictly negative: $K_{ij} < 0$ $\forall ij \in E$. Thus, $\dot{\hat{\phi}}_k = \sum_{l \in N_k} \hat{K}_{kl} \cdot f_{kl}(\hat{\phi}_l - \hat{\phi}_k) < 0$. Similarly, if $q_{ij} > q_{ij}^+$ $\forall ij \in E$, then $K_{ij} > 0$ $\forall ij \in E$, and $\dot{\hat{\phi}}_k = \sum_{l \in N_k} \hat{K}_{kl} \cdot f_{kl}(\hat{\phi}_l - \hat{\phi}_k) > 0$. In both cases $(\hat{\phi}, \hat{K})$ cannot be an equilibrium. $\square$

In [40] it was demonstrated that when in system (2) each $f_{ij} = \sin(\phi_j - \phi_i)$, $q_{ij} = q_{ij}^+ = 1$, and the underlying topology is a complete graph, then the in-phase equilibrium is asymptotically stable. Theorem 4 thus provides a complementary set of results to [40] under more general topologies and coupling functions.

A similar result can be formulated for the anti-phase equilibria of system (2).

**Theorem 5 (Stability of anti-phase equilibria)** *If functions $f_{ij}()$ satisfy Assumption 2, then an anti-phase equilibrium $(\phi^*, K^*)$ is:*

- **asymptotically stable**, *if $q_{ij} > q_{ij}^-$ (when $\phi_i^* = \phi_j^*$), and $q_{ij} < q_{ij}^+$ (when $|\phi_i^* - \phi_j^*| = \pi$), $\forall ij \in E$;*
- **unstable**, *if $q_{ij} < q_{ij}^-$ (when $\phi_i^* = \phi_j^*$), and $q_{ij} > q_{ij}^+$ (when $|\phi_i^* - \phi_j^*| = \pi$), $\forall ij \in E$.*

*Proof.* The proof is again based on the sufficient instability and stability conditions of Theorems 2 and 3. According to Theorem 3 (condition (27)), an anti-phase equilibrium $(\phi^*, K^*)$ is stable if

$$(F_{ij} + q_{ij}) \cdot f'_{ij} > 0, \forall ij \in E, \qquad (28)$$

because at such equilibria $f_{ij} = 0$ for any $i$ and $j$. If $\phi_i^* = \phi_j^*$, then $f'_{ij} > 0$ due to the Assumption 2. Therefore, condition (28) is satisfied if $q_{ij} > -F_{ij}(0)$ which is equivalent to $q_{ij} > -F_{ij}^{max} = q_{ij}^-$. Similarly, if $|\phi_i^* - \phi_j^*| = \pi$, then $f'_{ij} < 0$, and inequality (28) will be fulfilled if $F_{ij}(\pi) + q_{ij} < 0$, that is, if $q_{ij} < -F_{ij}^{min} = q_{ij}^+$.

From Theorem 2, an anti-phase equilibrium will be unstable if

$$(F_{ij} + q_{ij}) \cdot f'_{ij} < 0, \forall ij \in E, \qquad (29)$$

again since at such equilibria $f_{ij} = 0$ for any $i$ and $j$, and because $K_{ij} = \alpha_{ij}(F_{ij} + q_{ij})$. If $\phi_i^* = \phi_j^*$, then $f'_{ij} > 0$, and condition (29) is met if $F_{ij}(0) + q_{ij} < 0$, i.e. if $q_{ij} < -F_{ij}^{max} = q^-$. If in turn $|\phi_i^* - \phi_j^*| = \pi$, then $f'_{ij} < 0$, and equilibrium is unstable if $F_{ij}(\pi) + q_{ij} > 0$, i.e. if $q_{ij} > -F_{ij}^{min} = q_{ij}^+$. $\square$

Notice that the first part of Theorem 4 follows from Theorem 5.

**Remark** *Requirements for instability of Theorems 4 and 5 are generally more restrictive than the instability condition of Theorem 2. Indeed, these requirements guarantee that $K_{ij}f'_{ij} - \alpha_{ij}f_{ij}^2 < 0$ $\forall ij \in E$ which is a stronger requirement than (9).*

A direct consequence of Theorems 4 and 5 is the following result.

**Corollary 6** *If $q_{ij}$ takes values from an interval $(q_{ij}^-, q_{ij}^+)$ $\forall ij \in E$, then the in-phase and all anti-phase equilibria of system* (2) *are asymptotically stable.*

Conditions formulated in Theorems 2 and 3 are sufficient, and therefore, there may exist equilibria of system (2) whose stability cannot be characterized by these conditions. We overcome this problem in the next section by requiring underlying topology to be a tree graph. For such graphs we provide a criterion of stability that allows us to verify the stability of any equilibrium.

## 5 PHASE LOCKING FOR TREE TOPOLOGY

In this section we consider system (2) when the underlying topology graph $G$ is a tree. For example, star and chain graphs are two graphs belonging to this type of topology. We apply the results formulated in the previous section and show how they can be further extended for the tree graphs. We first consider in Subsection 5.1 the case of a general coupling, when parameters $q_{ij}$ are allowed to take any arbitrary values except for $q_{ij}^+$ and $q_{ij}^-$, and then explore a special case when each network connection is either strictly attractive ($q_{ij} > q_{ij}^+$) or strictly repulsive ($q_{ij} < q_{ij}^-$) in Subsection 5.2.

### 5.1 General Coupling

When the topology is a tree, each single edge defines a cut of $G$, and condition (9) for a single-edge cut $C(P) = ij$ can be written as $(F_{ij} + q_{ij}) \cdot f'_{ij} - f_{ij}^2 < 0$, because $K_{ij} = \alpha_{ij}(F_{ij} + q_{ij})$ at an equilibrium and $\alpha_{ij} > 0$. Thus, using Theorem 2, a sufficient instability condition for tree graphs can be formulated as follows.

**Corollary 7 (Sufficient instability condition for trees)** *If there exists an edge $ij \in E$ such that at equilibrium $(\phi^*, K^*)$ of system* (2) *with tree topology and functions $f_{ij}$ satisfying Assumption 1:*

$$(F_{ij} + q_{ij}) \cdot f'_{ij} - f_{ij}^2 < 0, \quad (30)$$

*where $F_{ij} = F_{ij}(\phi_j^* - \phi_i^*)$, $f'_{ij} = f'_{ij}(\phi_j^* - \phi_i^*)$ and $f_{ij} = f_{ij}(\phi_j^* - \phi_i^*)$, then $(\phi^*, K^*)$ is an unstable equilibrium.*

Using Theorem 3 and Corollary 7, the stability of an equilibrium of system (2) with a tree topology can be determined if $(F_{ij} + q_{ij}) \cdot f'_{ij} - f_{ij}^2 \neq 0$ for every $ij \in E$. As will be further shown, the last condition is always satisfied if all functions $f_{ij}()$ satisfy Assumption 2. Besides these additional assumptions on functions $f_{ij}()$, we also require throughout this subsection that $q_{ij} \neq q_{ij}^+ = -F_{ij}^{min}$ and $q_{ij} \neq q_{ij}^- = -F_{ij}^{max}$, $\forall ij \in E$. The following fact establishes a property of all equilibria of system (2) with a tree topology, and with functions $f_{ij}()$ satisfying Assumption 2.

**Lemma 8** *Let $(\phi^*, K^*)$ be an equilibrium of system* (2) *with a tree underlying topology and with functions $f_{ij}()$ satisfying Assumption 2, and let $f_{ij} = f_{ij}(\phi_j^* - \phi_i^*)$, $F_{ij} = F_{ij}(\phi_j^* - \phi_i^*)$, $K_{ij} = K_{ij}^*$. Suppose in addition, that $q_{ij} \neq q_{ij}^+$ and $q_{ij} \neq q_{ij}^-$, $\forall ij \in E$. Then, for each pair $ij \in E$ exactly one of the following two conditions is satisfied:*

- *$f_{ij} = 0$; this implies that either $\phi_j^* - \phi_i^* = 0$, or $\phi_j^* - \phi_i^* = \pi$.*
- *$K_{ij} = 0$; this implies that $F_{ij} + q_{ij} = 0$.*

*Proof.* Since the underlying topology is defined by a graph $G$ that is a tree, there are nodes in $G$ each of which has a single neighbor. These nodes are the leaves of a tree graph $G$. Let $\phi_i$ be an oscillator associated with a leaf node $i$, and let $\phi_j$ be an oscillator such that node $j$ is a single neighbor of node $i$. Then, from equation (2a), at an equilibrium $(\phi^*, K^*)$: $\dot{\phi}_i = 0$ if and only if $K_{ij} = 0$ or $f_{ij}(\phi_j^* - \phi_i^*) = 0$. Notice that $K_{ij} = 0$ if and only if $F_{ij} + q_{ij} = 0$. When $q_{ij} \neq -F_{ij}^{min}$ and $q_{ij} \neq -F_{ij}^{max}$, for function $f_{ij}$ satisfying Assumption 2, $f_{ij}$ and $F_{ij} + q_{ij}$ are not equal to zero simultaneously. This implies that either $f_{ij} = 0$ or $K_{ij} = F_{ij} + q_{ij} = 0$. We then remove all leaf nodes from the graph $G$ and apply the same reasoning for the leaves of a new smaller graph which is also a tree. We repeat this procedure until we obtain a single node, and at each step the condition of the theorem is satisfied. □

**Corollary 9** *Equilibria of system* (2) *with a tree topology and functions $f_{ij}()$ satisfying Assumption 2, are isolated.*

We can now see that when the underlying topology of system (2) is a tree, $q_{ij} \neq q_{ij}^+$, $q_{ij} \neq q_{ij}^-$ ($\forall ij \in E$), and if coupling functions $f_{ij}()$ satisfy Assumption 2, then at any equilibrium:

$$(F_{ij} + q_{ij}) \cdot f'_{ij} - f_{ij}^2 \neq 0 \quad (31)$$

for every $ij \in E$. Indeed, if at equilibrium $f_{ij} = 0$, then from Lemma 8, $K_{ij} \neq 0$, and $f'_{ij} \neq 0$ due to Assumption 2. Hence, $(F_{ij} + q_{ij}) \cdot f'_{ij} \neq 0$. If $K_{ij} = 0$, i.e. $F_{ij} + q_{ij} = 0$, then by definition of $F_{ij}$, properties of $q_{ij}$, and from Assumption 2: $f_{ij} \neq 0$.

Since Assumption 2 guarantees that condition (31) holds for every equilibrium, it is possible to formulate a criterion of stability for system (2) with underlying tree topology.

**Theorem 10 (Criterion of stability for tree topology)** *If the underlying topology of system* (2) *is a tree, $q_{ij} \neq q_{ij}^+$, $q_{ij} \neq q_{ij}^-$, and functions $f_{ij}()$ satisfy Assumption 2 for each $ij \in E$, then an equilibrium $(\phi^*, K^*)$ is stable if and only if condition* (27) *holds for every edge $ij \in E$. Moreover, each stable equilibrium is also asymptotically stable.*

*Proof.* Suppose that for equilibrium $(\phi^*, K^*)$ condition (27) is satisfied for any $ij \in E$, then this equilibrium is asymptotically stable by Theorem 3. Now assume that there exists an edge $ij \in E$ such that condition (27) does not hold for it, i.e.

$$(F_{ij} + q_{ij}) \cdot f'_{ij} - f_{ij}^2 \leq 0.$$

The above inequality must be strict since condition (31) holds for trees. Now consider a cut $C(P) = ij$ defined by a single edge $ij$. It immediately follows that the equilibrium must be unstable by Theorem 2. □

From this criterion we also conclude that (under the criterion's conditions) an equilibrium with $K_{ij}^* = 0$ for some edge $ij \in E$ will be unstable. We now provide a result regarding ranks of matrices $\hat{B}$ and $M$ for a tree topology. This result will be later used to show convergence to a stable equilibrium almost surely in Theorem 12.

**Lemma 11** *At any equilibrium of system* (2) *with a tree underlying topology and* $q_{ij} \neq q_{ij}^+$, $q_{ij} \neq q_{ij}^-$ ($\forall ij \in E$):

$$n = m + 1,$$
$$Rank(\hat{B}) = Rank(\hat{B}^T) = \min(n+m, m+m) = 2m,$$
$$Rank(M) = 2m.$$

*Proof.* First two equations are satisfied even outside of the equilibria. In particular, the first equation says that in a tree graph the number of vertices is greater than the number of edges by one, and is straightforward. The second fact follows from the properties of the incidence matrix $B$ and because we consider a tree topology. We now show that at an equilibrium the third equation is correct. Recall that matrix $M$ has dimensions $(2m) \times (2m)$. Suppose by contradiction that at some equilibrium $Rank(M) < 2m$, then there exists a non-zero vector $\vec{x}$ such that $M\vec{x} = 0$. This system contains $2m$ equations, the first $m$ of them are of the form:

$$x_i K_i f_i' + x_{m+i} f_i = 0, \quad (32)$$

where $i = 1, \ldots, m$, and the remaining $m$ equations are:

$$x_i \cdot f_i + x_{m+i}/\alpha_i = 0, \quad (33)$$

where $i = 1, \ldots, m$ again. We choose a particular index $i$ such that $1 \leq i \leq m$, and consider a pair of corresponding equations: one from (32) and another one from (33). According to Lemma 8, at an equilibrium exactly one of the following conditions holds: $f_i = 0$ or $K_i = 0$. We first consider the case when $f_i = 0$. Then, from the first equation, $x_i = 0$ because $K_i \neq 0$, $f_i' \neq 0$ and $f_i = 0$. And from the second equation: $x_{m+i} = 0$ since $\alpha_i > 0$. Now consider the case when $K_i = 0$, then $f_i \neq 0$. From the first equation: $x_{m+i} = 0$, and then from the second equation: $x_i = 0$ since $f_i \neq 0$. Therefore, in both cases $x_i = x_{m+i} = 0$, and since this should be true for all $1 \leq i \leq m$, vector $\vec{x}$ has to be a zero vector which contradicts our assumption that $\vec{x}$ is non-zero. □

While Theorem 1 does not guarantee pointwise convergence and isolation of equilibria in general, for the case of a tree underlying topology we proved isolation of equilibria in Corollary 9, and now can show that the system converges to a stable equilibrium almost surely.

**Theorem 12** *At any equilibrium of system* (2) *with a tree topology,* $q_{ij} \neq q_{ij}^+$, $q_{ij} \neq q_{ij}^-$ ($\forall ij \in E$), *and with functions $f_{ij}$ satisfying Assumption 2, the Jacobian has only one zero eigenvalue due to rotational invariance, and system converges to a stable equilibrium almost surely.*

*Proof.* The Jacobian matrix $J$ and symmetric matrix $A$ in decomposition (10) are both of size $(n+m) \times (n+m)$ or, using Lemma 11, of size $(2m+1) \times (2m+1)$. Moreover, matrix $A$ can be expressed as follows:

$$A = -\hat{B}M\hat{B}^T, \quad (34)$$

where matrix $M$ is of size $(2m) \times (2m)$, and matrices $\hat{B}$ and $\hat{B}^T$ are of size $(2m+1) \times (2m)$, $(2m) \times (2m+1)$, respectively. We now employ the following fact: if matrix $A_1$ has dimensions $x \times y$, matrix $A_2$ is of size $y \times z$ and $Rank(A_2) = y$, then $Rank(A_1 A_2) = Rank(A_1)$. Using this fact we conclude that $Rank(\hat{B}M) = 2m$, and $Rank(A) = Rank(\hat{B}M\hat{B}^T) = 2m$. Because rank of a symmetric matrix is equal to the number of its non-zero eigenvalues, $A$ has a single zero eigenvalue. Then, since $J$ and $A$ have the same numbers of zero eigenvalues (as was shown in the proof to Theorem 2), Jacobian matrix $J$ has only one zero eigenvalue which is due to the rotational invariance. Thus, all equilibria of (2) with a tree topology and functions $f_{ij}$ satisfying Assumption 2, are hyperbolic (i.e. they do not have any center manifold) when domain of the system is restricted to the subspace orthogonal to $[\mathbf{1}_n \ \mathbf{0}_m]^T$. Therefore, system (2) almost surely converges to a stable equilibrium (due to Proposition 1 of [12], for example). □

**Remark** *Theorem 12 implies that for a tree topology at least one stable equilibrium exists. We can verify this fact by constructing a stable in-phase or anti-phase equilibrium for a given system* (2) *with a tree underlying topology, functions $f_{ij}()$ satisfying Assumption 2, and parameters $q_{ij}$ such that $q_{ij} \neq q_{ij}^+$ and $q_{ij} \neq q_{ij}^-$ for all edges $ij \in E$. Without loss of generality, the phase value of each oscillator at such equilibrium is either 0 or $\pi$. To construct a stable equilibrium, we arbitrarily choose an oscillator to be the root of a tree, assign phase value equal to zero to this oscillator, and then traverse the tree using, for example, the breadth-first search. During the tree traversal, at each iteration we consider an edge connecting a previously visited oscillator with assigned phase value and a new oscillator whose phase is determined at this iteration. For example, if edge $ij \in E$ is considered, and previously visited oscillator $i$ has phase $\phi_i^* = 0$, then the phase of oscillator $j$ is assigned with value 0 if $q_{ij} > q_{ij}^+$, and with value $\pi$, otherwise. The equilibrium constructed in this manner will be asymptotically stable due to Theorem 5.*

### 5.2 Attractive and Repulsive Coupling

All previously formulated results also hold for systems where each connection $ij \in E$ is either strictly attractive ($q_{ij} > q_{ij}^+$) or strictly repulsive ($q_{ij} < q_{ij}^-$) coupling. Additionally, Lemma 8 can be further improved: for strictly repulsive or attractive connections, $K_{ij}$ cannot attain a zero value, and thus each equilibrium is characterized by conditions $f_{ij} = 0$ for each $ij \in E$. This implies that at each equilibrium all phase differences are multiples of $\pi$. Moreover, the following result holds.

**Theorem 13 (Almost global stability for trees)** *Suppose in system* (2) *the functions $f_{ij}()$ satisfy Assumption 2, and the underlying topology is a tree graph. If each connection $ij \in E$ is either strictly attractive ($q_{ij} > q_{ij}^+$) or strictly repulsive ($q_{ij} < q_{ij}^-$), then system* (2) *has a (unique) almost globally asymptotically stable in-phase or anti-phase equilibrium $(\phi^*, K^*)$ that, for each edge $ij \in E$, $(\phi^*, K^*)$ satisfies:*

$$\begin{cases} \phi_i^* = \phi_j^* & \text{if } q_{ij} > q_{ij}^+, \\ |\phi_i^* - \phi_j^*| = \pi & \text{if } q_{ij} < q_{ij}^-. \end{cases} \quad (35)$$

*Proof.* As was previously mentioned, under the conditions of this theorem the equilibrium set of system (2) consists of in-phase and anti-phase equilibria. Due to Theorem 10, an equilibrium $(\phi^*, K^*)$ is stable if and only if condition $(F_{ij} + q_{ij}) \cdot f_{ij}' > f_{ij}^2$ holds for each edge $ij \in E$. Because $f_{ij} = 0$ at each equilibrium of system (2) satisfying assumptions of Theorem 13, an equilibrium is stable if and only if $(F_{ij} + q_{ij}) \cdot f'_{ij} > 0$. If $q_{ij} > q^+_{ij}$, then this condition will be satisfied if and only if $f'_{ij} > 0$, which implies that $\phi^*_i = \phi^*_j$ due to Assumption 2. If $q_{ij} < q^-_{ij}$, then $f'_{ij}$ must be negative, and it means that $|\phi^*_i - \phi^*_j| = \pi$. Therefore, condition (35) is both necessary and sufficient for stability of an equilibrium. Moreover, according to Theorem 10, if an equilibrium is stable, it is also asymptotically stable.

Existence of a stable in-phase or anti-phase equilibrium satisfying (35) was shown by construction in the remark to Theorem 12. It remains to show, therefore, the uniqueness of this constructed stable equilibrium. Suppose, there exists another stable equilibrium $(\hat{\phi}, \hat{K})$, which is different from the constructed equilibrium $(\phi^*, K^*)$. Since it is stable, equilibrium $(\hat{\phi}, \hat{K})$, must also satisfy condition (35). And because $(\hat{\phi}, \hat{K})$ is different from $(\phi^*, K^*)$, there exists at least one pair of oscillators $k$ and $l$ (not necessary connected) such that $|\phi^*_k - \phi^*_l| \neq |\hat{\phi}_k - \hat{\phi}_l|$. It cannot happen if $k$ and $l$ are connected, i.e. if $kl \in E$, because of the condition (35). Now, assume that $k$ and $l$ are not connected, then since the underlying topology graph $G$ is a tree, there exists a *single* shortest path from $k$ to $l$ in $G$: $k, p_1, \ldots, p_h, l$, where $p_1, \ldots, p_h \in V$ are some oscillators and $h \geq 1$. This path consists of $h + 1$ edges $kp_1, p_1p_2, \ldots, p_hl$, and let $\mathcal{M}$ be the number of edges in this path with strictly repulsive coupling ($q < q^-$). Then, to satisfy condition (35), in any stable equilibrium the phase difference between oscillators $k$ and $l$ must be equal to zero if $\mathcal{M}$ is even, and equal to $\pi$ if $\mathcal{M}$ is odd. Thus, $|\phi^*_k - \phi^*_l|$ must be equal to $|\hat{\phi}_k - \hat{\phi}_l|$, which contradicts our assumption about equilibrium $(\hat{\phi}, \hat{K})$. □

**Corollary 14** *If for system* (2) *under the assumptions of Theorem 13, all connections are strictly attractive ($q_{ij} > q^+_{ij}$ $\forall ij \in E$), then the in-phase equilibrium is almost globally asymptotically stable. If all connections are strictly repulsive ($q_{ij} < q^-_{ij}$ $\forall ij \in E$), then the unique anti-phase equilibrium satisfying: $|\phi^*_i - \phi^*_j| = \pi$ $\forall ij \in E$ is almost globally asymptotically stable.*

Condition $|\phi^*_i - \phi^*_j| = \pi$ $\forall ij \in E$ implies that oscillators are divided into two sets (corresponding to phase values 0 and $\pi$, for example) so that each graph edge connects an oscillator from one set with an oscillators in the other set. This division can be done if the graph is bipartite which is always the case when the graph is a tree. Moreover, since the graph is connected, bipartition is unique, and thus, the anti-phase equilibrium satisfying condition $|\phi^*_i - \phi^*_j| = \pi$ $\forall ij \in E$ is also unique.

## 6 PHASE LOCKING FOR ARBITRARY TOPOLOGY

In Theorem 12 we demonstrated that when the underlying topology is a tree, convergence to a stable equilibrium occurs almost surely. It was possible to show this fact mainly due to the characterization of the equilibria formulated in Lemma 8. In particular, for the case of a tree topology all equilibria are isolated and moreover, hyperbolic. However, equilibria are not necessary isolated in the case of a non-tree topology as was demonstrated in Subsection 3.3 by means of an example of four completely connected oscillators. In that example all values of parameters $q_{ij}$ were equal to zero.

It turns out that the example in Subsection 3.3 is degenerate. That is, the non-isolation of the equilibria was mainly due to the special choice of $q_{ij}$ and $\alpha_{ij}$. Thus, we can formulate a result similar to Theorem 12 by making additional assumption on the choice of parameters $\alpha$ and $q$. This leads again towards a condition that guarantees convergence to a stable equilibrium almost surely.

**Theorem 15** *Suppose all functions $f_{ij}()$ satisfy Assumption 2, and for each $ij \in E$ parameters $\alpha_{ij}$ and $q_{ij}$ are chosen randomly from continuous probability distributions on the intervals $(0, \infty)$ and $(-\infty, \infty)$, respectively. Then, with probability one in the selection of these parameters, system* (2) *converges to a stable equilibrium almost surely.*

*Proof.* The proof is similar to the proof of Theorem 4 in [33] and is based on the parametric transversality theorem (Theorem 6.35 of [28]). To get rid of a zero eigenvalue of the Jacobian $J$ of system (2) which corresponds to rotational invariance, instead of $n$ phase variables $\phi_1, \ldots, \phi_n$ we will consider $n - 1$ phase differences $\mu_j \triangleq \phi_{j+1} - \phi_1$, where $j = 1, \ldots, n - 1$. Thus, all other phase values are measured relative to the phase value of the first oscillator. Notice, since the sum of all phases is an invariant of system (2), variables $\mu_j$, $j = 1, \ldots, n - 1$, uniquely define phase values $\phi_i$, where $i = 1, \ldots, n$. Let matrices $R \in \mathbb{R}^{n \times (n-1)}$ and $U \in \mathbb{R}^{n \times (n-1)}$ be defined as

$$R = \begin{bmatrix} \mathbf{0}^T_{n-1} \\ I_{n-1} \end{bmatrix}, U = \begin{bmatrix} -\mathbf{1}^T_{n-1} \\ I_{n-1} \end{bmatrix}, \quad (36)$$

where $I_{n-1}$ denotes an identity matrix of dimension $n - 1$, and $\mathbf{0}_{n-1}$, $\mathbf{1}_{n-1}$ are vectors of zeros and ones, respectively, of length $n - 1$. Note that $U^T R = I_{n-1}$ and $B^T = B^T R U^T$. Then, in the new variables $(\mu, K)$ system (2) can be rewritten as follows:

$$\begin{aligned} \dot{\mu} &= -U^T B \operatorname{diag}(K) f[B^T R \mu], \\ \dot{K} &= \mathcal{S}\Big(\mathcal{A}(F[B^T R \mu] + q) - K\Big). \end{aligned} \quad (37)$$

Here $F \in \mathbb{R}^m$ is a vector whose components are $F_{ij}$, and $q \in \mathbb{R}^m$ is a vector containing parameters $q_{ij}$ for each $ij \in E$. Further, $\mathcal{A} \in \mathbb{R}^{m \times m}$ and $\mathcal{S} \in \mathbb{R}^{m \times m}$ are diagonal matrices with parameters $\alpha_{ij}$ and $s_{ij}$ ($ij \in E$), respectively, on the diagonal. In the above equations and further in the proof, $\operatorname{diag}(x)$, where vector $x \in \mathbb{R}^m$, denotes a diagonal matrix of size $m \times m$ with elements of $x$ in its diagonal.

The main idea of the proof is to show that at the equilibria of system (37), the Jacobian does not have eigenvalues with zero real part. To show this we first prove that the Jacobian is invertible at each equilibrium for almost all values of parameters $\alpha$ and $q$, and then demonstrate that all its eigenvalues are real. Let $\mathcal{T}$ denote the finite collection of all $m \times m$ diagonal matrixes $\Delta = \operatorname{diag}\{\delta_1, \ldots, \delta_m\}$ such that either $\delta_k = (F_k(0) + q_k) \cdot f'_k(0)$ or $\delta_k = (F_k(\pi) + q_k) \cdot f'_k(\pi)$, where $k = 1, \ldots, m$ is the edge index. For each such matrix $\Delta \in \mathcal{T}$, we define the closed set

$$\mathcal{P}_\Delta = \{\alpha \in \mathbb{R}^m : \det(R^T B \mathcal{A} \Delta B^T R) = 0\}, \quad (38)$$

Matrix $\Delta$ is invertible due to Assumption 2 and because $(F_{ij}(x) + q_{ij}) \neq 0$ for $x = 0$ and $x = \pi$, since $q_{ij} \neq q^-_{ij}$ and $q_{ij} \neq q^+_{ij}$. Additionally, the columns of $B^T R$ are independent since they are rows $2, \ldots, n$ of the incidence matrix $B$. Then, $\mathcal{P}_\Delta \neq \mathbb{R}^m$ (for instance, if $\mathcal{A} = \Delta^{-1}$), and thus $\mathcal{P}_\Delta$ is a closed algebraic set having zero measure. Therefore, $\mathcal{P} = \bigcup_{\Delta \in \mathscr{T}} \mathcal{P}_\Delta$ is also a closed algebraic set having zero measure, and set $\mathcal{O}_\alpha = \mathbb{R}^m_+ \setminus \mathcal{P}$ is a nonempty open set. Further, let $\mathcal{O}_q$ be the set of vectors $q \in \mathbb{R}^m$ whose components satisfy $q_{ij} \neq q_{ij}^-$ and $q_{ij} \neq q_{ij}^+$. It is easy to show that the set $\mathbb{R}^m \setminus \mathcal{O}_q$ has a zero measure.

Let $\mathcal{H}(\mu, K, q)$ be a mapping $\mathbb{T}^{n-1} \times \mathbb{R}^m \times \mathcal{O}_q \to \mathbb{R}^{n-1+m}$:

$$\mathcal{H}(\mu, K, q) = \begin{bmatrix} -U^T B \operatorname{diag}(K) f[B^T R \mu] \\ \mathcal{S}\Big(\mathcal{A}(F[B^T R \mu] + q) - K\Big) \end{bmatrix}, \quad (39)$$

and $\alpha \in \mathcal{O}_\alpha$. Notice that $\mathcal{H}(\mu, K, q) = 0$ only at the equilibria of system (37). Next, the Jacobian of $\mathcal{H}(\mu, K, q)$ is

$$D\mathcal{H}(\mu, K, q) = \begin{bmatrix} \frac{\partial \mathcal{H}}{\partial \mu} & \frac{\partial \mathcal{H}}{\partial K} & \frac{\partial \mathcal{H}}{\partial q} \end{bmatrix}, \quad (40)$$

where

$$\frac{\partial \mathcal{H}}{\partial \mu} = \begin{bmatrix} -U^T B \operatorname{diag}(K * f'[B^T R \mu]) B^T R \\ -\mathcal{S}\mathcal{A} \operatorname{diag}(f[B^T R \mu]) B^T R \end{bmatrix}, \quad (41)$$

$$\frac{\partial \mathcal{H}}{\partial K} = \begin{bmatrix} -U^T B \operatorname{diag}(f[B^T R \mu]) \\ -\mathcal{S} \end{bmatrix}, \quad (42)$$

$$\frac{\partial \mathcal{H}}{\partial q} = \begin{bmatrix} 0_{(n-1)\times m} \\ \mathcal{S}\mathcal{A} \end{bmatrix}, \quad (43)$$

and $0_{x \times y}$ denotes a zero matrix of dimensions $x$ by $y$.

We now show that when $\mathcal{H}(\mu, K, q) = 0$, matrix $D\mathcal{H}(\mu, K, q)$ has full row rank for all $(\mu, K, \alpha, q) \in \mathbb{T}^{n-1} \times \mathbb{R}^m \times \mathcal{O}_\alpha \times \mathcal{O}_q$. Consider the following block matrix:

$$\mathcal{W} = \begin{bmatrix} I_{n-1} & 0_{(n-1)\times m} \\ L & 0_{m \times m} \\ \mathcal{A}^{-1} L & I_m \end{bmatrix}, \quad (44)$$

where

$$L = \mathcal{A} \operatorname{diag}(f[B^T R \mu])^+ \cdot \operatorname{diag}\Big((F[0_m] + q) * f'[0_m] \\ - (F[B^T R \mu] + q) * f'[B^T R \mu]\Big) B^T R, \quad (45)$$

here $(\cdot)^+$ denotes the Moore-Penrose pseudoinverse. Then, since $K = \mathcal{A}(F[B^T R \mu] + q)$ when $\mathcal{H}(\mu, K, q) = 0$,

$$D\mathcal{H}(\mu, K, q) \cdot \mathcal{W} = \\ \begin{bmatrix} -U^T B \mathcal{A} \Delta(\mu) B^T R & 0_{(n-1)\times m} \\ -\mathcal{S}\mathcal{A} \operatorname{diag}(f[B^T R \mu]) B^T R & \mathcal{S}\mathcal{A} \end{bmatrix}, \quad (46)$$

where $\Delta(\mu)$ is a diagonal matrix

$$\Delta(\mu) = \operatorname{diag}\Big((F[B^T R \mu] + q) * f'[B^T R \mu]\Big) \\ + \operatorname{diag}(f[B^T R \mu]) \cdot \operatorname{diag}(f[B^T R \mu])^+ \cdot \\ \cdot \operatorname{diag}\Big((F[0_m] + q) * f'[0_m] - (F[B^T R \mu] + q) * f'[B^T R \mu]\Big). \quad (47)$$

Because functions $f()$ satisfy Assumption 2, $f(x) = 0$ if and only if $x \in \{0, \pi\}$, and therefore, for any $\mu \in \mathbb{T}^{n-1}$, the matrix $\Delta(\mu)$ belongs to $\mathscr{T}$. It implies that matrix $R^T B \mathcal{A} \Delta(\mu) B^T R$ is invertible for $(\mu, K, \alpha, q) \in \mathbb{T}^{n-1} \times \mathbb{R}^m \times \mathcal{O}_\alpha \times \mathcal{O}_q$. Next, because

$$U^T B = U^T U R^T B, \quad (48)$$

and matrix $U^T U$ is invertible, it follows that matrix $-U^T B \mathcal{A} \Delta(\mu) B^T R = -(U^T U)(R^T B \mathcal{A} \Delta(\mu) B^T R)$ in (46) is also invertible. We now conclude that the whole matrix (46) is invertible, because its rows are independent ($\mathcal{S}\mathcal{A}$ is a diagonal matrix with positive numbers on the diagonal). As a consequence of this, when $\mathcal{H}(\mu, K, q) = 0$, the Jacobian matrix $D\mathcal{H}(\mu, K, q)$ must have independent rows for all $(\mu, K, \alpha, q) \in \mathbb{T}^{n-1} \times \mathbb{R}^m \times \mathcal{O}_\alpha \times \mathcal{O}_q$.

Thus, $\mathcal{H} \pitchfork \{0\}$, and it follows from the parametric transversality theorem ( [28], Theorem 6.35) that there exists a set $\mathcal{Y} \subset \mathcal{O}_q$ having zero measure such that if $q \in \mathcal{O}_q \setminus \mathcal{Y}$, then $\mathcal{H}_q \pitchfork \{0\}$, where $\mathcal{H}_q$ denotes the mapping $(\mu, K) \to \mathcal{H}(\mu, K, q)$. Therefore, $(n-1+m) \times (n-1+m)$ matrix

$$\begin{bmatrix} \frac{\partial \mathcal{H}}{\partial \mu} & \frac{\partial \mathcal{H}}{\partial K} \end{bmatrix} \quad (49)$$

is invertible when $\mathcal{H}(\mu, K, q) = 0$ for almost all values of parameters $\alpha$ and $q$.

Now observe that matrix (49) is the Jacobian of system (37). It remains to show that this matrix can have only real eigenvalues. Indeed, (49) is equal to a product of a diagonal positively definite matrix and a symmetric matrix:

$$\begin{bmatrix} \frac{\partial \mathcal{H}}{\partial \mu} & \frac{\partial \mathcal{H}}{\partial K} \end{bmatrix} = \begin{bmatrix} U^T U & 0_{(n-1)\times m} \\ 0_{m \times (n-1)} & \mathcal{S}\mathcal{A} \end{bmatrix} \begin{bmatrix} -A_1 & -A_2 \\ -A_3 & -\mathcal{A}^{-1} \end{bmatrix}, \quad (50)$$

where

$$A_1 = R^T B \operatorname{diag}(K * f'[B^T R \mu]) B^T R, \quad (51)$$

$$A_2 = R^T B \operatorname{diag}(f[B^T R \mu]), \quad (52)$$

and

$$A_3 = \operatorname{diag}(f[B^T R \mu]) B^T R. \quad (53)$$

We observe that the second matrix in the product (50) is indeed symmetric since $A_2^T = A_3$, and $A_1$ is symmetric.

To summarize the results, for almost all randomly selected system parameters $\alpha_{ij}$ and $q_{ij}$, the Jacobian of system (37) at each equilibrium cannot have eigenvalues with zero real part, and system (37) converges to a stable equilibrium almost surely, so does system (2). $\square$

**Remark** *Theorem 15 states that system* (2) *converges to a stable equilibrium almost surely for almost all values of parameters $\alpha$ and $q$. In other words, the set of parameters $\alpha$ and $q$ for which convergence to a stable equilibrium is not guaranteed, has zero measure. The theorem, however, does not provide a description of this zero-measure set and, therefore, it cannot be applied to guarantee convergence to a stable equilibrium for a given example with some fixed values of parameters. This is the main distinction between this theorem and its analog Theorem 12 for a tree topology: the latter result guarantees convergence for all values of the parameters (except for $q = q^+$ and $q = q^-$).*

In Corollary 14 we proved that the in-phase equilibrium is almost globally asymptotically stable in the case of strictly attractive coupling if the underlying topology is a tree. With additional requirements on the coupling functions $f_{ij}()$ it is possible to show almost globally asymptotically stability of the in-phase equilibrium for an arbitrary topology. In particular, we now assume that each function $f_{ij}$ satisfies

**Assumption 3** *There exists a parameter $b \in (0, \pi)$ such that functions $f_{ij}$ $\forall ij \in E$ satisfy:*

TABLE 1: Behavior of system (2) for various choices of its parameters. Parameters $s_{ij}$ are from $(0,\infty)$ in each case.

| Theorem | Functions $f_{ij}$ | Topology | Parameters $q_{ij}$ | Parameters $\alpha_{ij}$ | Synchronization/Convergence to |
|---|---|---|---|---|---|
| Theorem 1 | Assumption 1 | arbitrary | $(-\infty, \infty)$ | $(0, \infty)$ | frequency synchronization |
| Theorem 12 | Assumption 2 | tree | $q_{ij} \neq q_{ij}^+$, $q_{ij} \neq q_{ij}^-$ | $(0, \infty)$ | stable eq. a.s. |
| Theorem 13 | Assumption 2 | tree | $q_{ij} > q_{ij}^+$ or $q_{ij} < q_{ij}^-$ | $(0, \infty)$ | unique stable in/anti-phase eq. satisfying (35) a.s. |
| Theorem 15 | Assumption 2 | arbitrary | $\mathcal{O}_q \setminus \mathcal{Y}$ | $\mathcal{O}_\alpha$ | stable eq. a.s. |
| Theorem 16 | Assumption 3 | arbitrary | $\mathcal{O}_q \setminus \mathcal{Y}$, $q_{ij} > q_{ij}^+$ | $\mathcal{O}_\alpha$ | in-phase stable eq. a.s. |

1) Assumption 1;
2) $f'_{ij}(x) > 0$, $\forall x \in [0, b)$;
3) $f'_{ij}(x) < 0$, $\forall x \in (b, \pi]$.

Notice that Assumption 3 is stronger than Assumption 2 in a sense that if functions $f_{ij}$ satisfy Assumption 3, they also satisfy Assumption 2. Function whose graph is shown on the left side of Fig. 1 satisfies Assumption 3. The following theorem is a generalization of [30] to systems with plastic connectivity.

**Theorem 16 (Almost global stability)** *Suppose functions $f_{ij}$ satisfy Assumption 3 with $b \leq \frac{\pi}{n-1}$, and $q_{ij} > q_{ij}^+$ for each $ij \in E$, then with probability one in the selection of parameters $\alpha$ and $q$, the in-phase equilibrium of system (2) is almost globally asymptotically stable.*

*Proof.* The proof is almost identical to the proof of Theorem 6 in [30] and is not presented here due to the space limitations. The main idea is that if all functions $f_{ij}$ satisfy Assumption 3 and $q_{ij} > q_{ij}^+$ $\forall ij \in E$, then the in-phase equilibrium is the only stable equilibrium of system (2). Then the statement of the theorem follows from Theorem 15. □

**Remark** *The in-phase equilibrium will be also almost globally asymptotically stable if all functions $(-f_{ij})$ satisfy Assumption 3 with $b \leq \frac{\pi}{n-1}$ and $q_{ij} < q_{ij}^-$ $\forall ij \in E$.*

In Table 1 we summarized the convergence results of system (2) with arbitrary and tree underlying topology, and for various values of parameters $\alpha$ and $q$.

## 7 NUMERICAL ILLUSTRATIONS

In this section we consider several network examples and apply our results to explore their behavior and investigate stability of their equilibria. In these examples we will assume that $f_{ij} = \sin(\phi_j - \phi_i)$ and $F_{ij} = \cos(\phi_j - \phi_i)$ $\forall ij \in E$, then $q_{ij}^+ = 1$ and $q_{ij}^- = -1$ $\forall ij \in E$. Thus, system (2) becomes a generalized Kuramoto model with plastic coupling strengths. Notice that this choice of functions $f_{ij}$ guarantees that Assumptions 2 and 3 are satisfied. Additionally, *(strictly) positive* and *(strictly) negative* coupling corresponds to (strictly) attractive and (strictly) repulsive coupling, respectively.

### 7.1 Two Oscillators (Lemma 8, Theorems 10, 12, 13, 16)

Since two connected oscillators form a tree topology, we can apply the results of Lemma 8 and Theorems 10, 12. In particular, for a tree topology case we can easily find all equilibria of the system using our results (Lemma 8). We will assume that $\alpha_{12}$ takes an arbitrary positive value.

If $q_{12} = 0$, then Lemma 8 can be applied to describe all equilibria of system (2) with two oscillators. The first type of equilibria corresponds to condition $f_{12} = 0$ which implies that $\sin(\phi_2 - \phi_1) = 0$, and $K/\alpha_{12} = F_{12} + q_{12} = f'_{12} + q_{12} = \cos(\phi_2 - \phi_1) + q_{12} = \cos(\phi_2 - \phi_1)$. The second type of equilibria is defined by condition $\cos(\phi_2 - \phi_1) = 0$, i.e. $K = 0$, and $f_{12} = \sin(\phi_2 - \phi_1) = 1$. Stability of each equilibrium type can be verified using a criterion provided in Theorem 10. The criterion's stability condition (27) takes the following form:

$$(\cos(\phi_2 - \phi_1) + q_{12}) \cdot \cos(\phi_2 - \phi_1) > 0. \quad (54)$$

If $q_{12} = 0$, then (54) is satisfied for the first type of equilibria (in-phase (Fig. 2c-i) and anti-phase (Fig. 2c-ii)), making these equilibria stable. The second equilibrium type corresponding to $K = 0$ is unstable (Fig. 2c-iii) because (54) does not hold for it.

Next, if $q_{12} = 2$ or $q_{12} = -2$, then system (2) of two oscillators has only in-phase and anti-phase equilibria. We can use Theorem 13 to conclude that when $q_{12} = 2$, then the in-phase equilibrium is stable (Fig. 2a-i) and the anti-phase equilibrium is unstable (Fig. 2a-ii), and in the case when $q_{12} = -2$, the in-phase (Fig. 2b-i) and anti-phase (Fig. 2b-ii) equilibria are unstable and stable, respectively. We can observe that all results are in agreement with Fig. 2 as expected.

According to Theorem 12, system (2) converges to a stable equilibrium almost surely for each of the considered values of parameter $q_{12}$. Moreover, when $q_{12} = 2$ (strictly attractive coupling) or when $q_{12} = -2$ (strictly repulsive coupling) the system has a unique almost global stable in-phase (if $q_{12} = 2$) or anti-phase (when $q_{12} = -2$) equilibrium, as predicted by Theorem 13. In addition, since $\frac{\pi}{2} = b < \frac{\pi}{2-1} = \pi$, we can apply Theorem 16 when $q_{12} = 2$ to reestablish the almost global stability of the in-phase equilibrium.

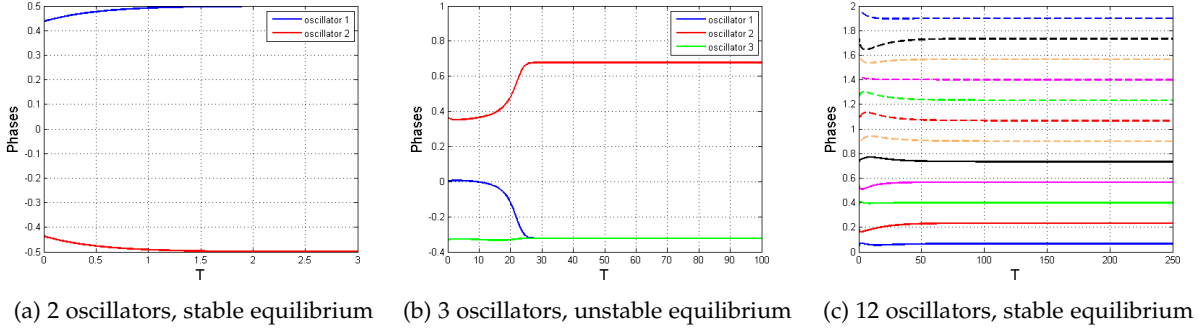

(a) 2 oscillators, stable equilibrium    (b) 3 oscillators, unstable equilibrium    (c) 12 oscillators, stable equilibrium

Fig. 5: Behavior of system (2) with $q_{ij} = 0$ ($\forall ij \in E$) of two (left), three (center) and twelve (right) oscillators after a small perturbation from a stable anti-phase, unstable and stable equilibrium, respectively.

The behavior of the system with $q_{12} = 0$ of two oscillators after a small perturbation from the anti-phase equilibrium is shown in Fig. 5a. This equilibrium is stable, and the system converges to it after a perturbation.

### 7.2 Three Oscillators (Theorems 2, 3, 4, 5, 15, 16)

Here we examine the stability of equilibria and the behavior of system (2) with three all-to-all connected oscillators. In this case the underlying topology is not a tree and we will employ Theorems 2 and 3 to show stability or instability. We will assume for simplicity that $\alpha_{ij} = 1$ for all $ij \in E$, and as in Subsection 3.2, will explore the cases $q_{ij} = 0$, $q_{ij} = 2$ and $q_{ij} = -2$ ($\forall ij \in E$).

We first consider the case when $q_{ij} = 0$ for each $ij \in E$. The in-phase equilibrium is stable (Fig. 3c-i) since condition (27) is satisfied: $\cos^2(0) > \sin^2(0)$. Clearly, the anti-phase equilibrium, i.e. when $\phi_1 = \phi_2$ and $\phi_3 = \phi_1 + \pi$, is also stable (Fig. 3c-ii). Now consider equilibrium $\phi_2 = \phi_1 + 2\pi/3$, $\phi_3 = \phi_1 - 2\pi/3$, and a two-edge cut $C(P) = \{12, 13\}$. Because $2(\cos^2(2\pi/3) - \sin^2(2\pi/3)) < 0$, condition (9) is satisfied and the equilibrium is unstable (Fig. 3c-iii). The next equilibrium is defined as $\phi_1 = \phi_2$, $\phi_3 = \phi_1 + \pi/2$. Using our cut $C(P) = \{13, 23\}$, we obtain: $2(\cos^2(\pi/2) - \sin^2(\pi/2)) < 0$, which means that the equilibrium is unstable (Fig. 3c-iv) due to Theorem 2. If the equilibrium is described by $\phi_2 = \phi_1 + \pi/3$, $\phi_3 = \phi_1 - \pi/3$, then using the cut $C(P) = \{12, 13\}$, we get $2(\cos^2(\pi/3) - \sin^2(\pi/3)) < 0$, and thus the equilibrium is unstable (Fig. 3c-v). In Fig. 5b behavior of this system is shown after a small perturbation from this unstable equilibrium. Finally, when $\phi_2 = \phi_1 + \pi/2$, $\phi_3 = \phi_1 - \pi/2$, we can use the cut $C(P) = \{12, 13\}$ to get $2(\cos^2(\pi/2) - \sin^2(\pi/2)) < 0$ which makes the equilibrium unstable (Fig. 3c-vi).

We now assume that $q_{ij} = 2, \forall ij \in E$. Then, the in-phase equilibrium is stable (Fig. 3a-i) since $(\cos(0)+2) \cdot \cos(0) > 0$, and condition (27) is satisfied. Next, the cut $C(P) = \{12, 13\}$ can be used to demonstrate the instability of the anti-phase equilibrium (Fig. 3a-ii): $2(\cos(\pi) + 2) \cdot \cos(\pi) < 0$, so that condition (9) is fulfilled. Using the same cut $C(P)$, the instability of equilibrium $\phi_2 = \phi_1 + 2\pi/3$, $\phi_3 = \phi_1 - 2\pi/3$ can be shown (Fig. 3a-iii).

Finally, we explore the case when $q_{ij} = -2$, $\forall ij \in E$. The in-phase equilibrium is unstable (Fig. 3d-i): $(\cos(0) - 2) \cdot \cos(0) < 0$. The stability of the anti-phase equilibrium, however, cannot be verified by our sufficient condition in Theorem 3, because $(\cos(0) - 2) \cdot \cos(0) < 0$. Further, equilibrium $\phi_2 = \phi_1 + 2\pi/3$, $\phi_3 = \phi_1 - 2\pi/3$ is stable (Fig. 3d-iii): $(\cos(2\pi/3) - 2) \cdot \cos(2\pi/3) - \sin^2(2\pi/3) > 0$. And Theorem 2 does not allow us to verify the instability of the last equilibrium (Fig. 3d-iv).

Stability or instability of in-phase and anti-phase equilibria in some cases can be also checked using Theorems 3 and 4. As was pointed out in the Remark of Subsection 4.3, the instability conditions of Theorems 4 and 5 are generally weaker than condition of Theorem 2. For example, while we could apply Theorem 2 to show that the anti-phase equilibrium is unstable when $q_{ij} = 2$ ($\forall ij \in E$), the sufficient instability condition of Theorem 5 is not fulfilled.

Because the topology of this example is not a tree, Theorem 12 cannot be applied. In addition, Theorem 15 formulated for an arbitrary topology, does not guarantee convergence to a stable equilibrium for given values of the parameters $\alpha$ and $q$. Nevertheless, when $q_{ij} = 2$, Theorem 16 can be used since $\frac{\pi}{2} = b = \frac{\pi}{3-1} = \frac{\pi}{2}$. Therefore, the in-phase equilibrium is almost globally asymptotically stable when $q_{ij} = 2$ ($\forall ij \in E$), which is in agreement with Fig. 3a.

### 7.3 Four Oscillators (Theorems 2, 15)

In Subsection 3.3 we described a set of equilibria of system (2) with $q_{ij} = 0$ ($\forall ij \in E$) and $\alpha_{ij} = \alpha > 0$ ($\forall ij \in E$) of four all-to-all connected oscillators, characterized by a parameter $\beta \in (0, \pi/2)$. We will first check that each equilibrium of this set is unstable. At any such equilibrium with fixed $\beta \in (0, \pi/2)$: $\phi_2 = \phi_1 + \pi/2$, $\phi_3 = \phi_1 + \beta$, $\phi_4 = \phi_1 - (\pi/2 - \beta)$, and we define a cut $C(P)$ as $C(P) = \{12, 13, 14\}$. Since $(\cos^2(\pi/2) - \sin^2(\pi/2)) + (\cos^2(\beta) - \sin^2(\beta)) + (\sin^2(\beta) - \cos^2(\beta)) < 0$, then Theorem 2 guarantees that this equilibrium is unstable.

This set of equilibria is of a special interest, because it consists of non-isolated equilibria, and the result of Theorem 15 does not apply to this system. Therefore, the values of parameters $\alpha$ and $q$ corresponding to this example, i.e. $\alpha_{ij} = \alpha > 0$ ($\forall ij \in E$) and $q_{ij} = 0$ ($\forall ij \in E$) belong to the zero-measure set of parameters that is excluded in the statement of Theorem 15. The requirement $\alpha_{ij} = \alpha > 0$ ($\forall ij \in E$) can be generalized: it is sufficient that $\alpha_{13} = \alpha_{14}$ and $\alpha_{23} = \alpha_{24}$ for the described set of equilibria to persist.

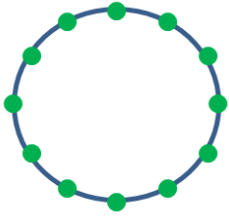

Fig. 6: Stable equilibrium corresponding to the example with twelve oscillators with ring topology.

### 7.4 Twelve Oscillators

In [40] it was shown that when in system (2) all coupling functions are $\sin()$, $q_{ij} = 0$ ($\forall ij \in E$), and the underlying topology is a complete graph, then the only stable equilibria are those in which every phase difference is a multiple of $\pi$. This property does not generally hold in the case of an arbitrary topology as we will demonstrate here using an example with twelve oscillators and $\sin()$ coupling functions. The underlying topology is a ring graph, so that the pairs of connected oscillators are $(1,2)$, $(2,3)$, $(3,4)$, ..., $(11,12)$, $(12,1)$. We will assume that $\alpha_{ij} = \alpha > 0$ ($\forall ij \in E$), and $q_{ij} = q$ ($\forall ij \in E$). The equilibrium is defined by the following phase values: $\phi_1 = 0$, and $\phi_i = \phi_{i-1} + \pi/12$, where $i = 2, \ldots, 12$. The phases of the oscillators at the equilibrium are shown in Fig. 6.

We can apply Theorem 3 to find stability conditions for this equilibrium under various values of $q$. Condition (27) in this example takes the form:

$$(\cos(\beta) + q) \cdot \cos(\beta) - \sin(\beta)^2 > 0, \quad (55)$$

where $\beta = \frac{2\pi}{n}$ is the phase difference between two neighboring oscillators in equilibrium, and $\beta = \frac{\pi}{6}$ when $n = 12$. The inequality (55) implies that if $q > -\frac{1}{\sqrt{3}}$, then the equilibrium is stable for $n = 12$ oscillators. In Fig. 5c the behavior of the system is shown after a small perturbation from this equilibrium.

## 8 Conclusion

In this work we studied a model of arbitrarily interconnected homogeneous coupled oscillators with a plastic coupling. We demonstrated that systems of oscillators described by this model always achieve frequency synchronization. Sufficient stability and instability conditions for equilibria were provided for a general underlying topology, and a criterion of stability was formulated for a tree topology. We additionally derived a sufficient condition that guarantees almost surely convergence to a stable equilibrium for tree topologies, and then obtained an analogous condition for the arbitrary topology case. Further, under certain assumptions on the coupling (strictly attractive or strictly repulsive connections), we formulated an almost global stability result for tree topologies. A similar condition was also derived for arbitrary topologies and strictly attractive connections. We illustrated our theoretical results with several examples.